\documentclass[11pt]{amsart}

\usepackage[utf8]{inputenc}
\usepackage{t1enc}
\usepackage{color}

\usepackage{amssymb}
\usepackage[all]{xy}

\usepackage{mathbbol}
\usepackage{latexsym}

\newcommand{\comments}[1]{}
\newcommand{\cst}{\textup{C}^*}
\newcommand{\id}{\operatorname{id}}
\newcommand{\I}{\operatorname{I}}
\newcommand{\Mor}{\operatorname{Mor}}
\newcommand{\Rep}{\operatorname{Rep}}

\newcommand{\flip}{\operatorname{flip}}
\newcommand{\Ad}{{\operatorname{Ad}}}

\newcommand{\tr}{\operatorname{tr}}

\newcommand{\M}{\operatorname{M}}

\newcommand{\tens}{\otimes}
\newcommand{\Tens}{\boxtimes}

\newcommand{\zero}{\{0\}}

\newcommand{\Ahat}{{\widehat A}}

\newcommand{\Ghat}{{\widehat G}}

\newcommand{\Rhat}{{\widehat R}}
\newcommand{\Vhat}{{\widehat V}}
\newcommand{\Deltahat}{{\widehat\Delta}}

\newcommand{\pihat}{\widehat{\pi}}

\newcommand{\vtilde}{{\widetilde v}}

\newcommand{\betatilde}{{\widetilde\beta}}

\newcommand{\Rtilde}{{\widetilde R}}
\newcommand{\alphatilde}{{\widetilde \alpha}}

\newcommand{\rtilde}{\widetilde{r}}

\newcommand{\zesp}{{\mathbb C}}

\newcommand{\real}{{\mathbb R}}

\newcommand{\inte}{{\mathbb Z}}

\newcommand{\Hil}{\mathcal{H}}

\newcommand{\skrC}{{\mathcal C}}
\newcommand{\skrL}{{\mathcal L}}

\newcommand{\skrD}{{\mathcal D}}
\newcommand{\skrR}{{\mathcal R}}

\newcommand{\skrU}{{\mathcal U}}

\newcommand{\Vs}[1]{\rule{0mm}{#1mm}}
\newcommand{\vs}{\vspace{2mm}}

\newcommand{\compcent}[1]{\vcenter{\hbox{$#1\circ$}}} 
\newcommand{\comp}{\mathbin{\mathchoice 
{\compcent\scriptstyle}{\compcent\scriptstyle} 
{\compcent\scriptscriptstyle}{\compcent\scriptscriptstyle}}}

\newcommand{\rf}[1]{{\rm (\ref{#1})}}

\newcommand{\set}[2]{\left\{#1:#2\right\}}

\newcommand{\xxx}[1]{}

\numberwithin{equation}{section}
\newtheorem{Thm}{Theorem}[section]
\newtheorem{Def}[Thm]{Definition}

\newtheorem{Prop}[Thm]{Proposition}

\newtheorem{rem}[Thm]{Remark}

\newenvironment{pf}[1][]{\vs
\noindent{\it Proof#1. }}{\qed}

\newcommand{\etyk}[1]{\label{#1}\stepcounter{equation}\tag{\theequation}}

\newcommand{\wiersz}[3]{
\begin{tabular}{c}#1\end{tabular}
&\begin{tabular}{c}#2\end{tabular}
&\begin{tabular}{c}#3\end{tabular}
\\ \hline}

\title{Monoidal category of $\cst$-algebras}
\author{S.L. Woronowicz}
\email{Stanislaw.Woronowicz@fuw.edu.pl}
\address{Instytut Matematyczny PAN
         ul.~\'Sniadeckich 8,
          00--956 Warszawa
          and
         Katedra Metod Mate\-matycznych Fizyki
         Wydzia\l Fizyki, Uniwersytet Warszawski,
         ul.~Pa\-steura 5,
         02-093 Warszawa,
         Poland.}

\date{June 2016}
\thanks{Supported by the Alexander von Humboldt-Stiftung. Partially supported by the National Science Center (NCN) grant no.~2015/17/B/ST1/00085.}

\begin{document}
\begin{abstract}
We consider the category of $\cst$-algebras equipped with actions of a locally compact quantum group. We show that this category admits a monoidal structure satisfying certain natural conditions if and only if the group is quasitriangular. The monoidal structures are in bijective correspondence with unitary $R$-matrices. To prove this result we use only very natural properties imposed on considered monoidal structures. We assume that monoidal product is a crossed product, monoidal product of injective morphisms is injective and that monoidal product reduces to the minimal tensor product when one of the involved $\cst$-algebras is equipped with a trivial action of the group. No a priori form of monoidal product is used.
\end{abstract}

\maketitle

\section{Introduction}
\label{Intro}

%Let $X$ and $Y$ be a norm closed subsets of a $\cst$ algebra. We set
%\[
%XY=\set{xy}{
%\begin{array}{c}
%x\in X\\y\in Y
%\end{array}}^{\rm CLS},
%\]
%

One of the paradigms of quantum theory says that the algebra of observables associated with a composed system is the tensor product of algebras associated with the parts of the system. This is the simplest monoidal functor. \vs

Given $\cst$-algebras $X$ and $Y$ one may consider the $\cst$-algebra $X\tens Y$. If $X$ and $Y$ are equipped with actions of a locally compact group $G$, then there exists unique action of $G$ on $X\tens Y$ such that natural embeddings of $X$ and $Y$ into $X\tens Y$ intertwine the actions of $G$. In the categorical language: tensor product $\tens$ defines a monoidal structure on the category $\cst_G$ of all $\cst$-algebras equipped with the action of $G$.\vs

This is no longer the case if $G$ is a quantum group. Let $G$ be a locally compact quantum group. We shall show that the category of $\cst_G$ admits a monoidal structure satisfying certain natural conditions if and only if the group is quasitriangular. The monoidal structures are in bijective correspondence with unitary $R$-matrices. In general the monoidal structure $\Tens$ does not coincide with $\tens$.\vs

Early examples of monoidal structures on $\cst_{G}$ are given in \cite{SLW95b} (for $G=\inte\times S^{1}$) and \cite{Row1} (for $G=\real$). To construct monoidal structure on $\cst_{G}$ was not a trivial task. Let $R$ be a unitary $R$-matrix. To define monoidal product $X\Tens Y$ of two objects $X,Y\in\cst_{G}$ one has to choose two representations $\alpha\in\Rep(X,K)$  and $\beta\in\Rep(Y,K)$ ($K$ is a Hilbert space) correlated in a way dictated by $R$. Then $X\Tens Y=\alpha(X)\beta(Y)$. The main problem consists in proving that $\alpha(X)\beta(Y)$ is a $\cst$-algebra i.e. that $\alpha(X)\beta(Y)=\beta(Y)\alpha(X)$. \vs

Methods used in \cite{Row1, SLW95b} took into account particular properties of the considered groups and gave no indication how to proceed in general case. A decisive step was made by Ryszard Nest and Christian Voigt in \cite{ NestVoight}. They showed that the intelligent use of Podle\'s condition (continuity of action) solves the problem. Nest and Voigt worked with locally compact quantum groups dual to Drinfeld doubles. The case of any quasitriangular locally compact quantum group was investigated in \cite{SLW16b}. The monoidal structure constructed in the latter paper have three very natural properties: monoidal product is a crossed product (Property 1), monoidal product of injective morphisms is injective (Property 2) and monoidal product reduces to the minimal tensor product when one of the involved $\cst$-algebras is equipped with a trivial action of the group (Property 3). See Proposition 4.6 in \cite{SLW16b}. \vs

In the present paper we show that any monoidal structure on $\cst_{G}$ with these properties is related to a unitary $R$-matrix in the way described in \cite{SLW16b}. What surprises in this result is the fact that it could be obtained in so general and abstract setting. We do not use any assumed in advance form of monoidal product. Instead we derive a compact formula relating monoidal structure with $R$-matrix. We show that the monoidal structure is uniquely determined by $R$-matrix.\vs

The basic notation used in the paper is recalled in section \ref{Sec1}. In particular we introduce a very general concept of crossed product of $\cst$-algebras. I turns out to be useful despite the fact that it admits degenerate cases which are very distant from the  original crossed product of an algebra by a group action.\vs

In our paper we consider $\cst$-algebras equipped with actions of a locally compact quantum group $G$. Although quantum groups plays a fundamental role in the subject, only the basic knowledge of the theory of locally compact quantum groups is required to understand this paper. All necessary informations are collected in section \ref{Sec2}. \vs

Section \ref{Sec3} contains the definition of a right action of a quantum group on a $\cst$-algebra. The notion of intertwining morphisms is introduced. These concepts lead to the category $\cst_{G}$. In section \ref{Sec4} we review categories, functors and natural mappings that appear in the paper. Some of them are defined in the section, the others, like $\Tens$ will be discussed later.\vs

In section \ref{Sec5} we investigate general properties of monoidal structures on $\cst_{G}$. Monoidal structure $\Tens$ is defined as a covariant associative functor acting from $\cst_{G}\times\cst_{G}$ into $\cst_{G}$. We assume that $\zesp$ is a neutral object for $\Tens$. The definition is followed by two natural mappings $\alpha$ and $\beta$ playing an important role in our considerations. In short $\alpha^{XY}$ and $\beta^{XY}$ denote natural embeddings of $X$ and $Y$ respectively into $X\Tens Y$. We establish a number of formulae involving $\alpha$ and $\beta$.
\vs
 
Monoidal structures considered in \cite{SLW95b,Row1,NestVoight,SLW16b} have some common very interesting properties. In section \ref{Sec6} we list these properties and investigate monoidal structures with these properties. Next we describe relation between monoidal structures and $R$-matrices. We formulate our main result and outline the proof. To prove that any  monoidal structure (obeying the Properties) is related to an $R$-matrix, one has to show that certain unitary element $\Rtilde\in\M(\Ahat\tens\Ahat\tens(A\Tens A))$ has trivial last leg (the one that corresponds to $A\Tens A$).

\vs 

To prove our main result we shall use two auxiliary propositions. They are discussed section \ref{Sec7}. Proposition \ref{cp} states that inside $X\Tens Y$, $G$-invariant elements of one algebra ($X$ or $Y$) commute with all elements of the other algebra. This fact is an easy consequence of Property 3. It turns out that the statement of Proposition \ref{cp} is equivalent to Property 3 (see Section \ref{Sec12}). In Proposition \ref{nogi} we deal with four $\cst$-algebras $X,Y,Z,T\in\cst_{G}$. Then $\M(X\Tens Z)$ and $\M(Y\Tens T)$ may be considered as subalgebras of $\M(X\Tens Y\Tens Z\Tens T)$. Proposition \ref{nogi} says that multiples of $\I_{X\Tens Y\tens Z\Tens T}$ are the only elements in the intersection $\M(X\Tens Z)\cap\M(Y\Tens T)$. In the case of tensor products (when `$\Tens$' is replaced by `$\tens$') the statement is obvious. It could be easily shown by using slicing maps $\omega\tens\id_{Y}\tens\mu\tens\id_{T}$ and $\id_{X}\tens\xi\tens\id_{Z}\tens\nu$ (where $\omega,\xi,\mu,\nu$ are continuous linear functionals on $X,Y,Z,T$ respectively). However for `$\Tens$' this technique (slicing maps) is not available and the proof of Proposition \ref{nogi} is more sophisticated.
\vs

Section \ref{Sec8} contains the proof of our main result. We show that the last leg of $\Rtilde$ satisfies the assumption of Proposition \ref{nogi}.\vs

Sections \ref{Sec9} -- \ref{Sec11} are devoted to the uniqueness of monoidal structure corresponding to a given $R$-matrix. To this end we investigate natural mappings $\Phi:\Tens\rightarrow\Tens'$, where $\Tens$ and $\Tens'$ are monoidal structures. Let $\id_{\cst_{G}}$ be identity functor acting on $\cst_{G}$. If $\varphi,\psi:\id_{\cst_{G}}\rightarrow\id_{\cst_{G}}$ are natural mappings then setting $\Xi^{XY}=\varphi^{X}\Tens\psi^{Y}$ we obtain a natural mapping from $\Tens$ into itself. Composition $\Psi\comp\Xi$ is another natural mapping from $\Tens$ into $\Tens'$. To avoid this ambiguity we introduce a concept of normalisation. \vs

We shall prove (cf Theorem \ref{j}) that any two monoidal structures $\Tens$ and $\Tens'$ on $\cst_{G}$ corresponding to the same $R$-matrix are related by unique normalized natural mapping $\Phi:\Tens\rightarrow\Tens'$. For any $X,Y\in\cst_{G}$, $\Phi^{XY}\in\Mor(X\Tens Y,X\Tens' Y)$ is an isomorphism. The existence of $\Phi$ means that the monoidal structure corresponding to a given $R$-matrix is unique.\vs

The proof of Theorem \ref{j} consists in two steps. First in section \ref{Sec10} we construct $\Phi^{AA}\in\Mor(A\Tens A, A\Tens' A)$, next (in section \ref{Sec11}) we extend the result to any pair $(X,Y)$ of objects of $\cst_{G}$. This extension is possible, because any object $X\in\cst_{G}$ is isomorphic to a subobject of $X\tens A$. The isomorphism is given by the action of $G$ on $X$. \vs

In section \ref{Sec12} we shall discuss alternative formulations of Property 3. It says that the monoidal product reduces to tensor product when one of involved $\cst$-algebras is equipped with trivial action of $G$. It turns out that one may restrict this demand to the situations when the product is taken with the two-dimensional $\cst$-algebra $D=\zesp^{2}$ equipped with trivial action of $G$.\vs

In section \ref{Sec13} we consider $\cst$-algebras equipped with left actions of $G$. We shall show, how to formulate our results in this context.

\section{Notation}
\label{Sec1}

Throughout the paper we shall use the following convention: if $T$ and $Z$ are norm closed subsets of a $\cst$ algebra then  $TZ$ will denote the closed linear span of the set of all products $tz$, where $t\in T$ and $z\in Z$:
\[
\etyk{conv}
TZ=\set{tz}{
\begin{array}{c}
t\in T\\z\in Z
\end{array}}^{\rm CLS},
\]
where CLS stands for norm Closed Linear Span.\vs

One of the basic categories considered in the paper is the category $\cst$ (see \cite{Vallin,SLW81}) whose objects are separable $\cst$ algebras. The morphisms are introduced in the following way: If $X,Y$ are $\cst$-algebras then  $\Mor(X,Y)$ is the set of all $^{*}$-algebra homomorphisms $\varphi$ acting from $X$ into $\M(Y)$ such that
$\varphi(X)Y=Y$. The latter formula uses the notation \rf{conv}.\vs

Any $\varphi\in\Mor(X,Y)$ admits a unique extension to a unital $^{*}$-algebra homomorphism acting from $\M(X)$ into $\M(Y)$. Composition of morphisms is defined as composition of their extensions. In what follows 
\[
\varphi:X\longrightarrow Y
\]
means that $\varphi\in\Mor(X,Y)$. It does not imply that $\varphi(X)\subset Y$.\vspace{5mm}

 In this paper for any $\cst$-algebras $X$ and $Y$, $X\tens Y$ always denote the minimal (spatial) tensor product. For any $x\in X$ and $y\in Y$ we set
 \[
\begin{array}{r@{\;=\;}l}
\alpha(x)&x\tens \I_{Y}\Vs{5}\\
\beta(y)&\I_{X}\tens y\Vs{5}
\end{array}
\]
Then $\alpha\in\!\Mor(X,X\tens Y)$, $\beta\in\!\Mor(Y,X\tens Y)$ and
 \[
 \alpha(X)\beta(Y)=X\tens Y.
 \]

 We shall use the following concept of crossed product algebra: Let $X,Y,Z$ be $\cst$-algebras, $\alpha\in\!\Mor(X,Z)$ and $\beta\in\!\Mor(Y,Z)$. We say that $Z$ is a crossed product of $X$ and $Y$ if
\[
\alpha(X)\beta(Y)=Z.
\]
\vs

%Example:

%{\bf Important:} In this paper for any $\cst$-algebras $X$ and $Y$, $X\tens Y$ always denote the minimal (spatial) tensor product.\vs

In practice crossed product of $\cst$-algebras appears in the way described in the following
\begin{Prop}
Let $X,Y$ be separable $\cst$ algebras, $\Hil$ be a Hilbert space, $\alpha\in\Rep(X,\Hil)$ and $\beta\in\Rep(Y,\Hil)$. Then
$
\alpha(X)\beta(Y)=\beta(Y)\alpha(X)
$
if and only if 
$
Z=\alpha(X)\beta(Y)$
 is a $\cst$-algebra.
Moreover in this case $\alpha\in\Mor(X,Z)$ and $\beta\in\Mor(Y,Z)$. Therefore $Z$ is a crossed product of $X$ and $Y$.
\end{Prop}

\begin{pf}
\[
\begin{array}{r@{\;=\;}l}
Z^{*}&(\alpha(X)\beta(Y))^{*}=\beta(Y)\alpha(X)=\alpha(X)\beta(Y)=Z,\\
ZZ&\alpha(X)\beta(Y)\alpha(X)\beta(Y)=\alpha(X)\alpha(X)\beta(Y)\beta(Y)\\
&\alpha(X)\beta(Y)=Z.
\end{array}
\]
It shows that $Z$ is $^{*}$-invariant and that $Z$ is closed with respect to the multiplication. Hence $Z$ is a $\cst$ algebra. Moreover we have 
\[
\begin{array}{c}
\alpha(X)Z=\alpha(X)\alpha(X)\beta(Y)=\alpha(X)\beta(Y)=Z,\\
\beta(Y)Z=\beta(Y)\beta(Y)\alpha(X)=\beta(Y)\alpha(X)=Z
\end{array}
\]
It shows that $\alpha\in\Mor(X,Z)$ and $\beta\in\Mor(Y,Z)$.
\end{pf}

\begin{Prop}
\label{P1}
Let 
\[
\etyk{D1}
\vcenter{
\xymatrix{
&Z\ar[d]&\\ 
X\ar[r]\ar[ru]^-{\alpha}\ar[rd]_-{\gamma}&S&Y\ar[l]\ar[lu]_-{\beta}\ar[ld]^-{\delta}\\
&T\ar[u]&
}}
\]
be a commutative diagram in the category $\cst$ such that vertical arrows are injective and $Z=\alpha(X)\beta(Y)$. Then there exists unique injective $\varphi\in\Mor(Z,T)$ such that
\end{Prop}

\[
\vcenter{
\xymatrix{
&Z\ar[dd]^-{\varphi}&\\ 
X\ar[ru]^-{\alpha}\ar[rd]_-{\gamma}&&Y\ar[lu]_-{\beta}\ar[ld]^-{\delta}\\
&T&
}}
\]
is a commutative diagram. If moreover $T=\gamma(X)\delta(Y)$ then $\varphi$ is an isomorphism.
 
\begin{pf}
Let $X',Y'\subset\M(S)$ be images of $X$ and $Y$ with respect to horizontal arrows and $T'=\gamma(X)\delta(Y)$. The diagram \rf{D1} shows that the images of $Z$ and $T'$ with respect to vertical arrows coincide with $X'Y'\subset\M(S)$. So they are equal, $T'$ is a $\cst$-algebra and composing (reading from right) $Z\rightarrow S$ with the inverse of $T'\rightarrow S$ we obtain the desired injection $\varphi: Z\rightarrow T$.
\end{pf}\vs

\section{Locally compact quantum groups}
\label{Sec2}
Let $G$ be a locally compact quantum group. This is a locally compact quantum space $G$ endowed with a continuous associative mapping $G\times G\longrightarrow G$ (group rule) subject to certain axioms.\vs

In practice we work with the $\cst$-algebra $A=\skrC_{0}(G)$ endowed with a morphism $\Delta\in\Mor(A,A\tens A)$ corresponding to the group rule on $G$. With shorthand notation:
\[
G=(A,\Delta).
\]
Strictly speaking one has to distinguish locally compact quantum group $G$ from the corresponding Hopf $\cst$-algebra $(A,\Delta)$. For instance
\[
\left\{
\begin{array}{c}
\text{actions}\\
\text{of}\ G
\end{array}
\right\}
=
\left\{
\begin{array}{c}
\text{coactions}\\
\text{of}\ (A,\Delta)
\end{array}
\right\}
\]
\vs

The present work does not use the full power of the Kustermans and Vaes theory (\cite{Vaes}, see also \cite{SLW03}) of locally compact quantum groups. Instead we use the theory of multiplicative unitaries (\cite{bs} and \cite{SLW96c}). For us locally compact quantum groups are objects coming from  manageable multiplicative unitary operators. In particular we do not use the Haar weights.\vs

Locally compact quantum groups appear in dual pairs:
\[
\begin{array}{r@{\;=\;}l}
G&(A,\Delta)\\
\Ghat&(\Ahat,\Deltahat)\Vs{5}
\end{array}
\]
The duality is described by a bicharacter $V$. This is a unitary element of $\M(\Ahat\tens A)$ such that 
\[
\etyk{birch}
\begin{array}{r@{\;=\;}l}
(\id\tens\Delta)V&V_{12}V_{13},\\
(\Deltahat\tens\id)V&V_{23}V_{13}.\Vs{5}
\end{array}
\]
One of the important feature of the theory of locally compact quantum groups is the unitary implementation of comultiplication: There exist faithful representations $\pi$ of $A$ and $\pihat$ of $\Ahat$ acting on the same Hilbert space $H$ such that for any $a\in A$ we have
\[
\etyk{Delta}
(\pi\tens\id_{A})\Delta(a)=V_{\pihat2}(\pi(a)\tens\I_{A})V_{\pihat2}^{*}.
\]
In this formula $V_{\pihat2}=(\pihat\tens\id)V$. We say that $(\pi,\pihat)$ is a Heisenberg pair.

%Let $(\pi,\pihat)$ be a Heisenberg pair for $G$. It means that $\pi$ and $\pihat$ are representations of $A$ and $\Ahat$ acting on the same Hilbert space $\Hil$ such that
%\[
%V_{\pihat3}V_{1\pi}=V_{1\pi}V_{13}V_{\pihat3}.
%\]
%Then for any $a\in A$ we have

\vs

Let $R$ be a unitary element of $\M(\Ahat\tens\Ahat)$.
We say that $R$ is a unitary $R$-matrix for $G$ if 
\[
\etyk{R}
\left\{\begin{array}{r@{\;=\;}l}
(\id\tens\Deltahat)R&R_{12}R_{13},\\
(\Deltahat\tens\id)R&R_{23}R_{13},\Vs{5}\\
R_{12}V_{13}V_{23}&V_{23}V_{13}R_{12}.\Vs{5}
\end{array}
\right.
\]

\begin{Def}
A locally compact quantum group $G=(A,\Delta)$ is called quasitriangular if there exists a unitary $R$-matrix in $\M(\Ahat\tens\Ahat)$.
\end{Def}

Let $R\in\M(\Ahat\tens\Ahat)$ be a unitary $R$-matrix and $\Rhat=\flip R^{*}$. Then $\Rhat\in\M(\Ahat\tens\Ahat)$ is a bicharacter. According to Theorem 5.3 of \cite{SLW12} there exists $\Delta_{R}\in\Mor(A,A\tens\Ahat)$ such that
\[
\etyk{DeltaR}
(\id_{\Ahat}\tens\Delta_{R})V=V_{12}\Rhat_{13}.
\]
\vs

We end this section with a short remark on opposite quantum groups. If $G=(A,\Delta)$ is a locally compact quantum group then, by definition $G^{\rm opp}=(A,\Delta^{\rm opp})$, where $\Delta^{\rm opp}=\flip\comp\Delta$. Consequently $\Ghat^{\rm opp}=(\Ahat,\Deltahat^{\rm opp})$, where $\Deltahat^{\rm opp}=\flip\comp\Deltahat$. Applying $^{*}$ to the both sides of \rf{birch} we get:
\[
\begin{array}{r@{\;=\;}l}
(\id\tens\Delta^{\rm opp})V^{*}&V_{12}^{*}V_{13}^{*},\\
(\Deltahat^{\rm opp}\tens\id)V^{*}&V_{23}^{*}V_{13}^{*}.\Vs{5}
\end{array}
\]
It shows that $\Ghat^{\rm opp}$ may be identified with the dual of $G^{\rm opp}$ with $V^{*}$ playing the role of bicharacter. If $G$ is quasitriangular and $R$ is the corresponding $R$-matrix then applying $^{*}$ to the both sides of \rf{R} we get:
\[
\begin{array}{r@{\;=\;}l}
(\id\tens\Deltahat^{\rm opp})R^{*}&R_{12}^{*}R_{13}^{*},\\
(\Deltahat^{\rm opp}\tens\id)R^{*}&R_{23}^{*}R_{13}^{*},\Vs{5}\\
R_{12}^{*}V_{13}^{*}V_{23}^{*}&V_{23}^{*}V_{13}^{*}R_{12}^{*}.\Vs{5}
\end{array}
\]
It shows that $G^{\rm opp}$ is quasitriangular with $R^{*}$ playing the role of $R$-matrix.

\section{$\cst$-algebras subject to an action of $G$}
\label{Sec3}
Let $X$ be a $\cst$-algebra and $\rho\in\Mor(X,X\tens A)$. We say that $\rho$ is an action of $G$ \linebreak on $X$ if

1. The diagram
\[
\etyk{a}
\vcenter{
\xymatrix{
X\ar[rr]^-{\rho}\ar[d]_-{\rho}&&X\tens A\ar[d]^-{\rho\tens\id}\\
X\tens A\ar[rr]_-{\id\tens\Delta}&&X\tens A\tens A
}}
\]
is a commutative,\vs

2. $\ker\rho=\zero$,\vs

3. $\rho(X)(\I\tens A)=X\tens A$ (Podle\'s condition).\vs

Condition 3 is a non-degeneracy condition of Podle\'s (cf \cite[Condition b of Definition 1.4]{Pod95}). For the first time Podle\'s condition appeared in his PhD dissertation\footnote{„Przestrzenie kwantowe i ich grupy symetrii”, available only in Polish.} in 1989. According to \cite{BSV}, Podle\'s condition characterises {\em strongly continuous} actions. \vs
 
{Remark:} Assume for the moment that $X$ and $A$ are algebras of operators acting on Hilbert spaces $K$ and $H$ respectively and that $\rho$ is a representation of $X$ acting on $K\tens H$. If $\rho(X)(\I\tens A)=X\tens A$ then $\rho(X)(X\tens A)=X\tens A$ and $\rho\in\Mor(X,X\tens A)$. In that sense
\[
\left(
\begin{array}{c}
\text{Podle\'s}\\ \text{condition}
\end{array}
\right)
\Longrightarrow
\left(\Vs{6}\rho\in\Mor(X,X\tens A)\right).
\]\vs

The main category considered in the paper is $\cst_{G}$. Objects of $\cst_{G}$ are $\cst$-algebras endowed with actions of $G$. For any $X\in\cst_{G}$, the action of $G$ on $X$ will be denoted by $\rho^{X}$. Morphisms in $\cst_{G}$ are $\cst$-morphisms intertwining the actions of $G$:
Let $X,Y\in\cst_{G}$. We say that a morphism $\gamma\in\Mor(X,Y)$ intertwins the actions of $G$ if
\[
\etyk{MorG}
\vcenter{
\xymatrix{
X\ar[rr]^{\gamma}\ar[d]_-{\rho^{X}}&&Y\ar[d]^-{\rho^{Y}}
\\
X\tens A\ar[rr]_-{\gamma\tens\id}&&Y\tens A
}}
\]
is a  commutative diagram. The set of all such morphisms will be denoted by $\Mor_{G}(X,Y)$. The reader should verify that for any $\varphi\in\Mor_{G}(X,Y)$ and $\psi\in\Mor_{G}(Y,Z)$ the composition $\psi\comp\varphi\in\Mor_{G}(X,Z)$ and that $\id_{X}\in\Mor_{G}(X,X)$.\vs

The following Proposition will be useful.\vs
 
\begin{Prop}
\label{4}
Let $X\in\cst_{G}$ and $u\in\M(X)$. Assume that $\rho^{X}(u)=\I_{X}\tens a$, where $a\in\M(A)$. Then $u=\lambda\I_{X}$ for some $\lambda\in\zesp$.
\end{Prop}
\begin{pf}
We have: $\I_{X}\tens\Delta(a)=(\id_{X}\tens\Delta)\rho^{X}(u)=$ $(\rho^{X}\tens\id_{A})\rho^{X}(u)=\I_{X}\tens\I_{A}\tens a$. Therefore $\Delta(a)=\I_{A}\tens a$ and (by known property of quantum groups) $a$ is a multiple of $\I_{A}$. Consequently $\rho^{X}(u)=\I_{X}\tens a$ is a multiple of $\I_{X\tens A}$. Remembering that $\rho^{X}$ is faithful we conclude that $u$ is a multiple of $\I_{X}$.
\end{pf}\vs

\section{Morphisms, functors and natural mappings}
\label{Sec4}

We shall use the language of
the theory of categories (see e.g. \cite{Semadeni}). Notions of
object, morphism, functor and natural mapping will appear. 
We work mainly with category $\cst_{G}$ introduced above.\vs

%The class of objects of the category will be denoted
%by the same symbol $\cst_{G}$ and the set of morphisms  acting from
%$X$ into $Y$ ($X,Y\in\cst_{G}$) will be denoted by
%$\Mor_{G}(X,Y)$.

We shall deal with the following functors and natural mappings:

\vspace{3mm}

\hspace*{1cm}
\begin{tabular}{|c|c|c|}\hline
\wiersz{  Functor  }{from}{  to  }
\wiersz{${\rm Proj}_1$\Vs{5}\\${\rm Proj}_2$\Vs{5}\\$\tens$\Vs{5}\\$\Tens$\Vs{5}\\$\Tens'$\Vs{5}}{  $\cst_{G}\times\cst_{G}$  }{  $\cst_{G}$  }
\wiersz{$\tens A$}{$\cst_{G}$\Vs{5}}{  $\cst_{G}$ \Vs{5} }
\end{tabular}
\hspace{10mm}
\begin{tabular}{|c|c|c|}\hline
\wiersz{Natural\\mapping}{from}{to}
\wiersz{$\alpha$}{${\rm Proj}_1$}{$\Tens$\Vs{5}}
\wiersz{$\beta$}{${\rm Proj}_2$}{$\Tens$\Vs{5}}
\wiersz{$\rho$}{$\id_{\cst_{G}}$\Vs{5}}{$\tens A$\Vs{5}}
\wiersz{$\Phi$}{$\Tens$\Vs{5}}{$\Tens'$\Vs{5}}
\end{tabular}\vspace{5mm}

%Objects of $\cst_{G}$ are $\cst$-algebras endowed with actions of $G$. For any $X\in\cst_{G}$, the action of $G$ on $X$ will be denoted by $\rho^{X}$. Morphisms in $\cst_{G}$ are $\cst$-morphisms intertwining the actions of $G$:
%Let $X,Y\in\cst_{G}$. We say that a morphism $\gamma\in\Mor(X,Y)$ intertwins the actions of $G$ if
%\[
%\etyk{MorG}
%\vcenter{
%\xymatrix{
%X\ar[rr]^{\gamma}\ar[d]_-{\rho^{X}}&&Y\ar[d]^-{\rho^{Y}}
%\\
%X\tens A\ar[rr]_-{\gamma\tens\id}&&Y\tens A
%}}
%\]
%is a  commutative diagram. The set of all such morphisms will be denoted by $\Mor_{G}(X,Y)$. The reader should verify that for any $\varphi\in\Mor_{G}(X,Y)$ and $\psi\in\Mor_{G}(Y,Z)$ the composition $\psi\comp\varphi\in\Mor_{G}(X,Z)$ and that $\id_{X}\in\Mor_{G}(X,X)$.\vs
% 

 {\bf Examples}

1.~Any $\cst$-algebra $X$ with the trivial action
\[
\rho^{X}(x)=x\tens \I_{A}\in \M(X\tens A)
\]
is an object of $\cst_{G}$.\vs

2.~The field of complex numbers $\zesp$ is a $\cst$-algebra. This is the initial object of category $\cst$: For any $\cst$-algebra $X$ the mapping
\[
1_{X}:\zesp\ni\lambda\longmapsto\lambda \I_{X}\in\M(X)
\]
is the only element of $\Mor(\zesp,X)$. Let $\rho^{\zesp}=1_{\zesp\tens A}$. Clearly $\rho^{\zesp}$ is a trivial action of $G$ on $\zesp$ and $\zesp\in\cst_{G}$.\vs

3.~The $\cst$-algebra $A=\skrC_{\infty}(G)$ with the action 
\[
\rho^{A}(a)=\Delta(a)\in \M(A\tens A)
\] 
is an object of $\cst_{G}$. This is a distinguished object.\vs

4.~Let $X$ be a $\cst$-algebra with any action of $G$. Then $X\tens A$ with the action
\[ 
\rho^{X\tens A}(x\tens a)=x\tens\Delta(a)\in \M((X\tens A)\tens A)
\]
is an object of $\cst_{G}$. The reader should notice that the action of $G$ on $X\tens A$ is induced by the action of $G$ on $A$. The action $\rho^{X}$ is ignored. However the commutative diagram \rf{a} shows that $\rho^{X}$ intertwines the actions of $G$ on $X$ and $X\tens A$:
\[
\rho^{X}\in\Mor_{G}(X,X\tens A).
\]

One may consider two functors: $\id_{\cst_{G}}$ and $\tens A$ (tensoring objects by $A$ and morphisms by $\id_{A}$) acting within $\cst_{G}$. Then $\rho$ become a natural mapping from  $\id_{\cst_{G}}$ into $\tens A$\vs

5.~Let $X,Y\in\cst_{G}$. Then $X\tens Y$ with the action
\[ 
\rho^{X\tens Y}(x\tens y)=x\tens\rho^{Y}(y)\in \M((X\tens Y)\tens A)
\]
is an object of $\cst_{G}$. Again the action $\rho^{X}$ is ignored: the action of $G$ on $X\tens Y$ is induced by the action of $G$ on $Y$. With the standard tensor product of morphisms, $\tens$ becomes a associative covariant functor acting from $\cst_{G}\times\cst_{G}$ into $\cst_{G}$.\vs

{\bf Remark}:
For any $X\in\cst_{G}$ we have three synonymous symbols
\[
\I_{X}\in\M(X), \hspace{4mm} \id_{X}\in\Mor_{G}(X,X),\hspace{4mm} 1_{X}\in\Mor_{G}(\zesp,X).
\]
\[
\etyk{zest}
 1_{\zesp}=\id_{\zesp}.
 \]
For any $\varphi\in\Mor_{G}(X,Y)$ we have
\[
\varphi(\I_{X})=\I_{Y}, \hspace{4mm} \varphi\comp\id_{X}=\varphi,\hspace{4mm} \varphi\comp1_{X}=1_{Y}.
\]

\section{Monoidal structures}
\label{Sec5}
\begin{Def}
%\label{}
A monoidal structure on the category $\cst_{G}$ is an associative covariant functor $\Tens$ acting from $\cst_{G}\times\cst_{G}$ into $\cst_{G}$ having $\zesp$ as neutral object.
\end{Def}

Being covariant functor means that $\Tens$ is a binary operation defined on objects and morphisms of $\cst_{G}$. For any $X,Y\in\cst_{G}$ we have an object $X\Tens Y\in\cst_{G}$. Moreover for any $X,Y,X',Y'\in\cst_{G}$ and any $\varphi\in\Mor_{G}(X,X')$ and $\psi\in\Mor_{G}(Y,Y')$ we have a morphism $\varphi\Tens\psi\in\Mor_{G}(X\Tens Y,X'\Tens Y')$. Composition of morphisms is compatible with $\Tens$:
\[
(\varphi'\Tens\psi')\comp(\varphi\Tens\psi)=(\varphi'\comp\varphi)\Tens(\psi'\comp\psi)
\]
for any $\varphi\in\Mor_{G}(X,X')$, $\varphi'\in\Mor_{G}(X',X'')$, $\psi\in\Mor_{G}(Y,Y')$ and $\psi'\in\Mor_{G}(Y',Y'')$ (where $X,X',X'',Y,Y',Y''\in\cst_{G}$).\vs

Associativity means that for any $X,Y,Z,X',Y',Z'\in\cst_{G}$ and any $\varphi\in\Mor_{G}(X,X')$, $\psi\in\Mor_{G}(Y,Y')$ and $\chi\in\Mor_{G}(Z,Z')$ we have:
\[
\etyk{asso}
\left\{
\begin{array}{r@{(}c@{\,\Tens\,}c@{)\,\Tens\,}c@{\;=\;}c@{\,\Tens\,(}c@{\,\Tens\,}c@{)}l}
&X&Y&Z&X&Y&Z&\\
&X'&Y'&Z'&X'&Y'&Z'&\\
&\varphi&\psi&\chi&\varphi&\psi&\chi&
\end{array}
\right.
\]
In what follows we shall omit brackets.\vs

``$\zesp$ is a neutral object'' means that
\[
\etyk{neutral}
\left\{
\begin{array}{c@{\;\Tens\;}c@{\;=\;}c@{\;=\;}c@{\;\Tens\;}c} 
X&\zesp&X&\zesp& X,\\
X'&\zesp&X'&\zesp& X',\\
\varphi&\id_{\zesp}&\varphi&\id_{\zesp}&\varphi
\end{array}
\right.
\]
for any $X,X'\in\cst_{G}$ and $\varphi\in\Mor(X,X')$.\vs

Except the case when $G$ is  the one-element group, the associative functor $\tens$ (see example 5 in the previous section) does not define a monoidal structure on $\cst_{G}$. This is because $\zesp$ is not a neutral object for $\tens$. Indeed, for any $X\in\cst_{G}$ we have: 
\[
\begin{array}{r@{\;=\;}l}
\zesp\tens X&X\\
X\tens\zesp&X_{tr},
\end{array}
\]
where $X_{\tr}$ is the $\cst$-algebra $X$ equipped with the trivial action of $G$. If $G$ is not trivial then $A_{\tr}\neq A$.\vs

The main result of this paper states that the 
category $\cst_{G}$ admits a monoidal structure (with certain natural properties) if and only if $G$ is quasi-triangular\footnote{the "if" part of the statement was established in \cite{SLW16b}}. More than that: monoidal structures are in one to one correspondence with unitary $R$-matrices.\vs

% {Natural mappings $\alpha$ and $\beta$}
Let $\Tens$ be a monoidal structure on $\cst_{G}$. For any $X,Y\in\cst_{G}$ we set
\[
\etyk{dfg}
\begin{array}{r@{\;=\;}c@{\:\Tens\:}c}
\alpha^{XY}&\id_{X}&1_{Y},\\
\beta^{XY}&1_{X}&\id_{Y}.\Vs{4}
\end{array}
\]
Then
\[
\begin{array}{r@{\;\in\Mor_{G}(}c@{,X\Tens Y)}c}
\alpha^{XY}&X&,\\
\beta^{XY}&Y&.
\end{array}
\]
In particular $\alpha^{X\zesp}\in\Mor(X,X)$ and $\beta^{\zesp Y}\in\Mor(Y,Y)$. Clearly
\[
\begin{array}{r@{\;=\;}l}
\alpha^{X\zesp}&\id_{X},\\
\beta^{\zesp Y}&\id_{Y}.
\end{array}
\]
Indeed using \rf{zest} and \rf{neutral} we have $\alpha^{X\zesp}=\id_{X}\Tens 1_{\zesp}=\id_{X}\Tens\id_{\zesp}=\id_{X}$ and similarly $\beta^{\zesp Y}=1_{\zesp}\Tens\id_{Y}=\id_{\zesp}\Tens\id_{Y}=\id_{Y}$.\vs

Let $X,Y,Z\in\cst_{G}$. Inserting in the equality
\[
(\varphi\Tens\psi)\Tens\chi=\varphi\Tens(\psi\Tens\chi)
\]
$\varphi$ equal either $\id_{X}$ or $1_{X}$, $\psi$ equal either $\id_{Y}$ or $1_{Y}$ and $\chi$ equal either $\id_{Z}$ or $1_{Z}$ and using \rf{dfg} we obtain six interesting equalities involving morphisms $\alpha$ and $\beta$:
\[
\begin{array}{r@{\;=\;}l@{\hspace{15mm}}r@{\;=\;}l}
\alpha^{X\Tens Y,Z}&\id_{X}\Tens\alpha^{YZ},&\beta^{X\Tens Y,Z}&\beta^{X,Y\Tens Z}\comp\beta^{YZ},\\
\Vs{5}\alpha^{XY}\Tens\id_{Z}&\id_{X}\Tens\beta^{YZ},&\alpha^{X\Tens Y,Z}\comp\beta^{XY}&\beta^{X,Y\Tens Z}
\comp\alpha^{YZ},\\
\Vs{5}\beta^{XY}\Tens\id_{Z}&\beta^{X,Y\tens Z},&\alpha^{X\Tens Y,Z}\comp\alpha^{XY}&\alpha^{X,Y\Tens Z}.
\end{array}
\] 

\begin{Prop}
%\label{}
Let $X,Y\in\cst_{G}$. Then the diagram
\[
\etyk{Natrho}
\vcenter{
\xymatrix{
X\ar[rrr]^-{\alpha^{XY}}\ar[dd]_-{\rho^{X}}&&&X\Tens Y\ar[dd]^-{\rho^{X\Tens Y}}&&&Y\ar[dd]^-{\rho^{Y}}\ar[lll]_-{\beta^{XY}}
\\
&&&&&&
\\
X\tens A\ar[rrr]_-{\alpha^{XY}\tens\,\id_{A}}&&&(X\Tens Y)\tens A&&&Y\tens A\ar[lll]^-{\beta^{XY}\tens\,\id_{A}}
}}
\]
is commutative.

\end{Prop}

\begin{pf}
Inserting in \rf{MorG}, $X\Tens Y$ instead of $Y$ and then setting $\varphi$ equal to $\alpha^{XY}$ and next equal to $\beta^{XY}$ we obtain \rf{Natrho}
\end{pf}

\begin{Prop}
%\label{}
$\alpha$ and $\beta$ are natural mappings from ${\rm Proj}_1$ and ${\rm Proj}_2$ into $\Tens$. More explicitly, for any $X,X',Y,Y'\in\cst_{G}$, $\varphi\in\Mor_{G}(X,X')$ and $\psi\in\Mor_{G}(Y,Y')$ the
diagram
\[
\etyk{Nat}
\vcenter{
\xymatrix{
X\ar[rr]^-{\alpha^{XY}}\ar[dd]_{\varphi}&&X\Tens Y\ar[dd]^{\varphi\Tens\psi}&&Y\ar[dd]^{\psi}\ar[ll]_-{\beta^{XY}}
\\
&&&&
\\
X'\ar[rr]_-{\alpha^{X'Y'}}&&X'\Tens Y'&&Y'\ar[ll]^-{\beta^{X'Y'}}
}}
\]
is commutative.
\end{Prop} 
\begin{pf}
We have: \newline 
\Vs{6}\indent $(\varphi\Tens\psi)\comp\alpha^{XY}=(\varphi\Tens\psi)\comp(\id_{X}\Tens\,1_{Y})=\varphi\Tens 1_{Y'}=(\id_{X'}\Tens\,1_{Y'})\comp\varphi=\alpha^{X'Y'}\comp\varphi$, \hfil and\newline 
\indent $(\varphi\Tens\psi)\comp\beta^{XY}=(\varphi\Tens\psi)\comp(1_{X}\Tens\id_{Y})=1_{X'}\Tens\psi=(1_{X'}\Tens\id_{Y'})\comp\psi=\beta^{X'Y'}\comp\psi$.

\end{pf}

\section{Natural properties of a monoidal structure on $\cst_{G}$}
\label{Sec6}

Let $G=(A,\Delta)$ be a quasitriangular locally compact quantum group. In \cite{SLW16b}  we introduced a monoidal structure $\Tens$ on the category $\cst_{G}$. 
It has the following properties:\vs

{\bf Property 1:}
For any $X,Y\in\cst_{G}$, $X\Tens Y$ is a crossed product of $X$ and $Y$:
\[
\etyk{wtorek}
X\Tens Y=\alpha^{XY}(X)\beta^{XY}(Y)
\]
 
{\bf Property 2:}
The $\Tens$-product of injective morphisms is injective.\vs
 
{\bf Property 3:}
$\Tens$ reduces to $\tens$, when the action of $G$ on one of the involved $\cst$-algebras is trivial. More precisely: 
If $X,Y\in\cst_{G}$ and if one of the actions $\rho^{X}$ and $\rho^{Y}$ is trivial then $X\Tens Y=X\tens Y$ as $\cst$-algebras. Moreover in this case
\[
\etyk{albet}
\begin{array}{r@{\;=\;}l}
\alpha^{XY}(x)&x\tens\I_{Y},\\
\beta^{XY}(y)&\I_{X}\tens y
\end{array}
\]
for any $x\in X$ and $y\in Y$.\vs

In what follows we shall be interested only in monoidal structures obeying Properties 1, 2 and 3.\vs

Let $X,Y\in\cst_{G}$. If $\rho^{X}$ is trivial then $X\Tens Y=X\tens Y$, $\alpha^{XY}$ and $\beta^{XY}$ are of the form \rf{albet} and diagram \rf{Natrho} shows that 
\[
\etyk{0911}
\begin{array}{l}
\rho^{X\Tens Y}(x\tens y)=x_{1}\rho^{Y}(y)_{23}\in X\tens Y\tens A,\\
\rho^{Y\Tens X}(y\tens x)=\rho^{Y}(y)_{13}x_{2}\in Y\tens X\tens A.\Vs{5}
\end{array}
\]
The above formulae hold for any $x\in X$ and $y\in Y$.
\vs

Let $X,Y,X',Y'\in\cst_{G}$, $\varphi\in\Mor_{G}(X,X')$ and $\psi\in\Mor_{G}(Y,Y')$. Assume that in each pair $(\rho^{X},\rho^{Y})$ and $(\rho^{X'},\rho^{Y'})$ one of the action is trivial. Then (by Property 3) $X\Tens Y=X\tens Y$, $X'\Tens Y'=X'\tens Y'$, morphisms $\alpha^{XY}$, $\alpha^{X'Y'}$, $\beta^{XY}$ , $\beta^{X'Y'}$ are of the form \rf{albet} and diagram \rf{Nat} shows that
\[
\etyk{trmor}
\varphi\Tens\psi=\varphi\tens\psi.
\]

{\bf Flip isomorphism:} 
Let $X,Y\in\cst_{G}$. Assume that one of the action $\rho^{X}$ and $\rho^{Y}$ is trivial. Then (by Property 3)
\[
\begin{array}{r@{\;=\;}l}
X\Tens Y&X\tens Y,\\
Y\Tens X&Y\tens X,
\end{array}
\]
as $\cst$-algebras (the corresponding actions of $G$ does not coincide).
In this case the flip map: $\flip^{XY}(x\tens y)=y\tens x$ may be considered as mapping acting from $X\Tens Y$ on $Y\Tens X$. Formulae \rf{0911} show that $\flip^{XY}$ intertwines the actions $\rho^{X\Tens Y}$ and $\rho^{Y\Tens X}$:
\[
\flip^{XY}\in\Mor_{G}(X\Tens Y,Y\Tens X).
\]
\vs

{\bf Mixed products.} Let $X,Y,Z\in\cst_{G}$ and $X_{\tr}$ be $X$ equipped with the trivial action of $G$:
$X_{\tr}=X\tens\zesp$. For trivial action $\Tens$ reduces to $\tens$. Therefore\vspace{-3mm}

\[
X\tens Y=X_{\tr}\Tens Y
\]
and by associativity
\[
\begin{array}{r@{\;=\;}l}
(X\tens Y)\Tens Z&(X_{\tr}\Tens Y)\Tens Z
\\&X_{\tr}\Tens (Y\Tens Z)=X\tens (Y\Tens Z).
\end{array}
\]
This way we showed the equality of mixed products
\[
\etyk{mixed}
(X\tens Y)\Tens Z=X\tens(Y\Tens Z).
\]
Clearly the similar formula holds for morphisms:
If $\varphi\in\Mor(X,X')$, $\phi\in\Mor_{G}(Y,Y')$ and $\psi\in\Mor_{G}(Z,Z')$ then
\[
\etyk{mix}
\left\{
\begin{array}{r@{\;=\;}l}
(X\tens Y)\Tens Z&X\tens(Y\Tens Z),\\
(X'\tens Y')\Tens Z'&X'\tens(Y'\Tens Z'),\\
(\varphi\tens\phi)\Tens\psi&\varphi\tens(\phi\Tens\psi).
\end{array}
\right.
\]
In particular
\[
\etyk{mi}
\left\{
\begin{array}{r@{\;=\;}c@{\tens\,}l}
\alpha^{X\tens Y,Z}&\id_{X}&\alpha^{YZ},\\
\beta^{X\tens Y,Z}&1_{X}&\beta^{YZ}.
\end{array}
\right.
\]
To get the first relation we put (in \rf{mix}) $\varphi=\id_{X}$, $\phi=\id_{Y}$ and  $\psi=1_{Z}$. Inserting $\varphi=1_{X}$, $\phi=1_{Y}$ and  $\psi=\id_{Z}$ we obtain the second relation. \vs

%\section{Monoidal structure and $R$-matrix}

The category $\cst_{G}$ contains a distinguished object $A$ with $\rho^{A}=\Delta$. Let $V\in\M(\Ahat\tens A)$ be the bicharacter describing the duality between $\Ghat$ and $G$. To make our formulae simpler we shall use the following shorthand notation:
\[
\etyk{abr}
\begin{array}{r@{\;=\;}l}
V_{1\alpha}&\left[\left(\id\tens\alpha^{AA}\right)V\right]_{13},\\
\Vs{6}V_{2\beta}&\left[\left(\id\tens\beta^{AA}\right)V\right]_{23}.
\end{array}
\]

Clearly $V_{1\alpha},V_{2\beta}\in \M(\Ahat\tens\Ahat\tens (A\Tens A))$. With this notation we have\footnote{see Proposition 4.6 in \cite{SLW16b}}:

\begin{Thm}
\label{MRW}
Let $G=(A,\Delta)$ be a quasitriangular locally compact quantum group with a unitary $R$-matrix $R\in \M(\Ahat\tens\Ahat)$. Then there exists  a monoidal structure $\Tens$ on $\cst_{G}$ having Properties 1, 2 and 3 and such that
\[
\etyk{RH}
V_{1\alpha}V_{2\beta}=V_{2\beta}V_{1\alpha}R_{12}.
\]
\end{Thm}

We shall prove the following:\vs

\begin{Thm}
\label{S}
Let $G=(A,\Delta)$ be a locally compact quantum group and $\Tens$ be a monoidal structure on $\cst_{G}$ having Properties 1, 2 and 3. Then there exists (unique) unitary $R$-matrix $R\in\M(\Ahat\tens\Ahat)$ such that
\[
V_{1\alpha}V_{2\beta}=V_{2\beta}V_{1\alpha}R_{12}.
\]
\end{Thm}

 {\bf Plan of the proof of Thm \ref{S}:}
Let 
\[
\Rtilde=V_{1\alpha}^{*}V_{2\beta}^{*}V_{1\alpha}V_{2\beta}.
\]
Then $\Rtilde\in\M(\Ahat\tens\Ahat\tens(A\Tens A))$. To prove Thm \ref{S} we have to show that the $(A\Tens A)$ - leg of $\Rtilde$ is trivial. In other words we have to show that $\Rtilde=R_{12}$, where $R\in\M(\Ahat\tens\Ahat)$. Next we have to prove that $R$ satisfies the relations \rf{R} characteristic for unitary $R$-matrix.

\section{Auxiliary statements}
\label{Sec7}

 Let $X\in\cst_{G}$ and $x\in\M(X)$. We say that $x$ is $G$-invariant if $\rho^{X}(x)=x\tens\I_{A}$. The proof of Theorem \ref{S} is based on the following two propositions: \vs

\begin{Prop}
\label{cp}
Let $X,Y\in\cst_{G}$, $x\in\M(X)$ and $y\in\M(Y)$. Assume that one of the elements $x$, $y$ is $G$-invariant.
Then 
\[
\alpha^{XY}(x)\beta^{XY}(y)=\beta^{XY}(y)\alpha^{XY}(x).
\]
\end{Prop}

\begin{Prop}
\label{nogi}
Let $X,Y,Z,T\in\cst_{G}$ and $u\in\M(X\Tens Z)$ and $v\in\M(Y\Tens T)$. Assume that
\[
\left(\id_{X}\Tens\;1_{Y}\Tens\id_{Z}\Tens\;1_{T}\right)(u)=
\left(1_{X}\Tens\id_{Y}\Tens\;1_{Z}\Tens\id_{T}\right)(v).
\]
Then $u=\lambda \I_{X\Tens Z}$ and $v=\lambda \I_{Y\Tens T}$, where $\lambda\in\zesp$.
\end{Prop}

\begin{pf}[ of Proposition \ref{cp}]
Assume for the moment that one of the actions $\rho^{X}$ and $\rho^{Y}$ is trivial. Then $\Tens$ becomes $\tens$: $X\Tens Y=X\tens Y$, 
\[
\begin{array}{r@{\;=\;}l}
\alpha^{XY}(x)&x\tens I_{Y},\\
\beta^{XY}(y)&I_{X}\tens y
\end{array}
\]
and $\alpha^{XY}(x)$ and $\beta^{XY}(y)$ obviously commute.\vs

Assume now that $x$ is a $G$-invariant element of $M(X)$. Let $X'$ be the smallest $\cst$-subalgebra of $M(X)$ containing $x$ and $I_{X}$. We provide $X'$ with the trivial action of $G$. Then $X'\in\cst_{G}$ and the embedding:
\[
\iota:X'\hookrightarrow\M(X)
\]
is a morphism in $\cst_{G}$: $\iota\in\Mor_{G}(X',X)$. Inserting in \rf{Nat} $\varphi=\iota$ and $\psi=\id_{Y}$ we see that
\[
\begin{array}{r@{\;=\;}l}
(\iota\Tens\id_{Y})\alpha^{X'Y}(x)&\alpha^{XY}(x),\\
(\iota\Tens\id_{Y})\beta^{X'Y}(y)&\beta^{XY}(y).
\end{array}
\]
By the first part of the proof $\alpha^{X'Y}(x)$ and $\beta^{X'Y}(y)$ commute. The above formulae show that $\alpha^{XY}(x)$ and $\beta^{XY}(y)$ commute.
The case, when $y$ is $G$-invariant may be treated in the same way.
\end{pf}\vs

\begin{pf}[ of Proposition \ref{nogi}]
Let
\[
\begin{array}{r@{\;=\;}l}
\varphi&\id_{X}\Tens\;1_{Y}\Tens\id_{Z}\Tens\;1_{T},\\
\psi&1_{X}\Tens\id_{Y}\Tens\;1_{Z}\Tens\id_{T}.
\end{array}
\]
Then
\[
\begin{array}{r@{\;\in\Mor_{G}(}c@{,X\Tens Y\Tens Z\Tens T)}l}
\varphi&X\Tens Z&,\\
\psi&Y\Tens T&.
\end{array}
\]
We assumed that
\[
\etyk{1}
\varphi(u)=\psi(v).
\]
At first we shall prove that
\[
\etyk{2}
(\rho^{Y}\Tens\id_{T})(v)=I_{Y}\tens\vtilde,
\]
where $\vtilde\in\M(A\Tens T)$.\vs

Let $\chi\in\Mor_{G}(Y,Y')$. Applying $\id_{X}\Tens\;\chi\Tens\id_{Z}\Tens\id_{T}$ to the both sides of \rf{1} we get:\hspace{5mm} $\left(\id_{X}\Tens\;1_{Y'}\Tens\id_{Z}\Tens\;1_{T}\right)(u)=
\left(1_{X}\Tens\id_{Y'}\Tens\;1_{Z}\Tens\id_{T}\right)(\chi\Tens\id_{T})(v)$. It shows that $(\chi\Tens\id_{T})(v)$ is independent of $\chi$, it depends only on $Y'$ - the target object  of $\chi$.\vs

Now, take faithful $\pi\in\Mor(Y,B_{0}(\Hil))$ (where $\Hil$ is a Hilbert space and $B_{0}(\Hil)$ is the algebra of all compact  operators on $\Hil$) and set $\chi=(\tau\pi\tens\id_{A})\rho^{Y}$, where $\tau$ is an automorphism of $B_{0}(\Hil)$. Then $\chi\in\Mor_{G}(Y,B_{0}(\Hil)\tens A)$. By the previous remark $(\chi\Tens\id_{T})(v)=(\tau\pi\tens\id_{A\Tens T})(\rho^{Y}\Tens\id_{T})(v)$ does not depend on $\tau$. Remembering that multiple of $\I_{\Hil}$ are the only operators invariant under all $\tau$ we obtain \rf{2}.\vs
 
We know that any morphism in the category $\cst_{G}$ intertwines the actions of $G$. In particular for 
\[
\begin{array}{r@{\,\in\Mor_{G}(}c@{,X\Tens Y\Tens Z)}l}
\id_{X}\Tens 1_{Y}\Tens\id_{Z}&X\Tens Z&,\\
1_{X}\Tens\id_{Y}\Tens 1_{Z}&Y&, \Vs{5}
\end{array}
\]
we have
\[
\begin{array}{r@{\rho^{X\Tens Y\Tens Z}(}c@{)=[(}c@{)\tens\id}c@{]\rho}l}
&\id_{X}\Tens 1_{Y}\Tens\id_{Z}&\id_{X}\Tens 1_{Y}\Tens\id_{Z}&_{A}&^{X\Tens Z},\\
&1_{X}\Tens\id_{Y}\Tens 1_{Z}&1_{X}\Tens\id_{Y}\Tens 1_{Z}&_{A}&^{Y}.\Vs{5}
\end{array}
\]
Tensoring ($\Tens$) from the right by $1_{T}$ the morphisms appearing in the first formula and by $\id_{T}$ the morphisms appearing in the second formula we get
\[
\begin{array}{r@{(\rho^{X\Tens Y\Tens Z}\Tens\id_{T})}c@{\,=\,}l}
&\varphi&\left[(\id_{X}\Tens 1_{Y}\Tens\id_{Z})\tens\alpha^{AT}\right]\rho^{X\Tens Z},            \\
&\psi&\left[(1_{X}\Tens\id_{Y}\Tens 1_{Z})\tens\id_{A\Tens T}\right](\rho^{Y}\Tens\id_{T}).\Vs{5}
\end{array}
\]
Applying the morphisms appearing in the first formula to $u$ and in the second formula to $v$ and using \rf{1} and \rf{2} we obtain:

\[
\etyk{3}
\left[\Vs{4}(\id_{X}\Tens 1_{Y}\Tens\id_{Z})\tens\alpha^{AT}\right]\rho^{X\Tens Z}(u)=\I_{X\Tens Y\Tens Z}\tens\vtilde.
\]
For $\tens$-product the technique of slices is available. Let $\omega$ be a state on $ X\Tens Y\Tens Z$. Then $\omega'=\omega\comp(\id_{X}\Tens 1_{Y}\Tens\id_{Z})$ is a state on $X\Tens Z$. Applying $\omega\tens\id_{A\Tens T}$ to the both sides of \rf{3} we see that
\[
\vtilde=\alpha^{AT}(a),
\]
where $a=\left(\Vs{3.5}\omega'\tens\id_{A}\right)\rho^{X\Tens Z}(u)\in\M(A)$. Comparing now the obvious formula
\[
\left[\Vs{4}(\id_{X}\Tens 1_{Y}\Tens\id_{Z})\tens\alpha^{AT}\right](\I_{X\tens Z}\tens a)=\I_{X\Tens Y\Tens Z}\tens\vtilde.
\]
with \rf{3} we conclude that $\rho^{X\Tens Z}(u)=\I_{X\Tens Z}\tens a$. Proposition \ref{4} shows now that $u$ is a multiple of $\I_{X\Tens Z}$. Proposition \ref{nogi} is proven.
\end{pf}

\section{Proof of Theorem \ref{S}}
 \label{Sec8}
\begin{pf}
Let $\Rtilde$ be a unitary element of $\M(\Ahat\tens\Ahat\tens(A\Tens A))$ introduced by the formula
\[
\Rtilde=V_{1\alpha}^{*}V_{2\beta}^{*}V_{1\alpha}V_{2\beta}.
\]
We have to show that the $(A\Tens A)$-leg of $\Rtilde$ is trivial. We shall deal with the $\Tens$-products of two and four copies of the distinguished object $A$:
\[
\begin{array}{c@{\;=\;}c@{\,\Tens\,}c}
A^{\Tens2}&A&A\\
A^{\Tens4}&A^{\Tens2}&A^{\Tens2}
\end{array}
\]
To make our formulae shorter we shall write $\alpha$ and $\beta$ instead of $\alpha^{AA}$ and $\beta^{AA}$ (this notation is coherent with \rf{abr}) and $\alphatilde$ and $\betatilde$ instead of $\alpha^{A^{\Tens2}A^{\Tens2}}$ and $\beta^{A^{\Tens2}A^{\Tens2}}$.
Then
\[
\begin{array}{c@{,\,}c@{\,\in\Mor_{G}(}c@{,\,}c@{)}c}
\alpha&\beta&A&A^{\Tens2}&,\\
\alphatilde&\betatilde&A^{\Tens2}&A^{\Tens4}&.
\end{array}
\]
Composing these morphisms we obtain four morphisms from $A$ into $A^{\Tens4}$. Using the formulae expressing $\alpha$ and $\beta$ as $\Tens$-products of $\id_{A}$ and $1_{A}$ one can easily verify that 
\[
\begin{array}{c@{\;=\;}c@{,\hspace{10mm}}c@{\,=\,}c}
(\alpha\Tens\alpha)\alpha&\alphatilde\alpha&(\beta\Tens\beta)\alpha&\alphatilde\beta,\\
(\alpha\Tens\alpha)\beta&\betatilde\alpha&(\beta\Tens\beta)\beta&\betatilde\beta.
\end{array}
\]

The following eight unitaries belonging to $\M(\Ahat\tens\Ahat\tens A^{\Tens4})$ will be involved in our computations:
For $i\in\{1,2\}$ and $r,s\in\{\alpha,\beta\}$ we set:
\[
V_{i,\rtilde s}=\left\{(\id_{\Ahat}\tens \rtilde\comp s)V\Vs{4}\right\}_{i3}.
\]
With this notation
\[
\etyk{7}
\begin{array}{r@{\;=\;}l}
\left(\id_{\Ahat\tens\Ahat}\tens(\alpha\Tens\alpha)\right)\Rtilde&V_{1,\alphatilde\alpha}^{*}V_{2,\betatilde\alpha}^{*}V_{1,\alphatilde\alpha}V_{2,\betatilde\alpha},\\
\left(\id_{\Ahat\tens\Ahat}\tens(\beta\Tens\beta)\right)\Rtilde&V_{1,\alphatilde\beta}^{*}V_{2,\betatilde\beta}^{*}V_{1,\alphatilde\beta}V_{2,\betatilde\beta}.
\end{array}
\]

We compute:
\[
\begin{array}{r@{\;=\;}l}
\left(\id_{\Ahat}\tens\rho^{A\Tens A}\right)V_{1\alpha}&\left(\id_{\Ahat}\tens\alpha\tens\id_{A}\right)(\id_{\Ahat}\tens\rho^{A})V\\\Vs{5}&\left(\id_{\Ahat}\tens\alpha\tens\id_{A}\right)V_{12}V_{13}=V_{1\alpha}V_{13}.
\end{array}
\]
Similarly
\[
\left(\id_{\Ahat}\tens\rho^{A\Tens A}\right)V_{1\beta}=V_{1\beta}V_{13}.
\]
Therefore
\[
\left(\id_{\Ahat}\tens\rho^{A\Tens A}\right)\left(V_{1\alpha}V_{1\beta}^{*}\right)=V_{1\alpha}V_{1\beta}^{*}\tens\I_{A}.
\]
It shows that the `second leg' of $V_{1\alpha}V_{1\beta}^{*}$ is $G$-invariant.\vs

Proposition \ref{cp} shows now that $V_{2,\betatilde\alpha}V_{2,\betatilde\beta}^{*}$ commutes with $V_{1,\alphatilde\beta}$ and that $V_{1,\alphatilde\alpha}V_{1,\alphatilde\beta}^{*}$ commutes with $V_{2,\betatilde\alpha}$. Using this information we get:
\[
\begin{array}{r@{\;=\;}l}
V_{1,\alphatilde\beta}^{*}V_{2,\betatilde\beta}^{*}V_{1,\alphatilde\beta}V_{2,\betatilde\beta}&
V_{1,\alphatilde\beta}^{*}V_{2,\betatilde\beta}^{*}V_{1,\alphatilde\beta}V_{2,\betatilde\beta}
V_{2,\betatilde\alpha}^{*}V_{2,\betatilde\alpha}=
V_{1,\alphatilde\beta}^{*}V_{2,\betatilde\alpha}^{*}V_{1,\alphatilde\beta}V_{2,\betatilde\alpha}\\
\Vs{5}&V_{1,\alphatilde\alpha}^{*}V_{1,\alphatilde\alpha}V_{1,\alphatilde\beta}^{*}V_{2,\betatilde\alpha}^{*}V_{1,\alphatilde\beta}V_{2,\betatilde\alpha}=V_{1,\alphatilde\alpha}^{*}V_{2,\betatilde\alpha}^{*}V_{1,\alphatilde\alpha}V_{2,\betatilde\alpha},
\end{array}
\]
We showed that the unitaries appearing on the right hand side of relations \rf{7} are equal. Therefore
\[
\left(\id_{\Ahat\tens\Ahat}\tens(\alpha\Tens\alpha)\right)\Rtilde=\left(\id_{\Ahat\tens\Ahat}\tens(\beta\Tens\beta)\right)\Rtilde.
\]
Notice that
\[
\begin{array}{r@{\;=\;}l}
\alpha\Tens\alpha&\id_{A}\Tens\, 1_{A}\Tens\id_{A}\Tens\, 1_{A},\\
\beta\Tens\beta&1_{A}\Tens\id_{A}\Tens\, 1_{A}\Tens\id_{A}.
\end{array}
\]
Proposition \ref{nogi} shows now that the `last leg' of $\Rtilde$ is trivial: $\Rtilde=R_{12}$, where\linebreak $R\in\M(\Ahat\tens\Ahat)$. To end the proof we have to show that $R$ satisfies \rf{R}. We already know that 
\[
\etyk{5}
V_{1\alpha}V_{2\beta}=V_{2\beta}V_{1\alpha}R_{12},
\]
Applying $\Deltahat$ to the first and second leg we get:
\[
\begin{array}{r@{\;=\;}l}
V_{2\alpha}V_{1\alpha}V_{3\beta}&V_{3\beta}V_{2\alpha}V_{1\alpha}\left\{(\Deltahat\tens\id_{\Ahat})R\right\}_{123},\\
V_{1\alpha}V_{3\beta}V_{2\beta}&V_{3\beta}V_{2\beta}V_{1\alpha}\left\{(\id_{\Ahat}\tens\Deltahat)R\right\}_{123}.\Vs{7}
\end{array}
\]
On the other hand we have:
\[
\begin{array}{r@{\;=\;}l}
V_{2\alpha}V_{1\alpha}V_{3\beta}&V_{2\alpha}V_{3\beta}V_{1\alpha}R_{13}\\
&V_{3\beta}V_{2\alpha}R_{23}V_{1\alpha}R_{13}=V_{3\beta}V_{2\alpha}V_{1\alpha}R_{23}R_{13},\\
V_{1\alpha}V_{3\beta}V_{2\beta}&V_{3\beta}V_{1\alpha}R_{13}V_{2\beta}\\
&V_{3\beta}V_{1\alpha}V_{2\beta}R_{13}=V_{3\beta}V_{2\beta}V_{1\alpha}R_{12}R_{13}
\end{array}
\]

It shows that
\[
\begin{array}{r@{\;=\;}l}
(\Deltahat\tens\id_{\Ahat})R&R_{23}R_{13},\\
(\id_{\Ahat}\tens\Deltahat)R&R_{12}R_{13}.     \Vs{6}
\end{array}
\]

To prove the third relation of \rf{R} we apply $\id_{\Ahat}\tens\id_{\Ahat}\tens\rho^{A\Tens A}$ to the both sides of \rf{5}:
\[
\begin{array}{c}
\left\{(\id_{\Ahat}\tens\rho^{A})V\right\}_{1\alpha3}\left\{(\id_{\Ahat}\tens\rho^{A})V\right\}_{2\beta3}\hspace*{20mm}\\ \hspace*{20mm}=
\left\{(\id_{\Ahat}\tens\rho^{A})V\right\}_{2\beta3}\left\{(\id_{\Ahat}\tens\rho^{A})V\right\}_{1\alpha3}R_{12}\Vs{6},
\end{array}
\]
\[
V_{1\alpha}V_{13}V_{2\beta}V_{23}=V_{2\beta}V_{23}V_{1\alpha}V_{13}R_{12},
\]
\[
V_{1\alpha}V_{2\beta}V_{13}V_{23}=V_{2\beta}V_{1\alpha}V_{23}V_{13}R_{12},
\]
\[
R_{12}V_{13}V_{23}=V_{23}V_{13}R_{12},\hspace*{5mm}
\]
\end{pf}

\section{Uniqueness of monoidal structures}
\label{Sec9}
Let $\Tens,\Tens'$ be monoidal structures on $\cst_{G}$ and $\Phi:\Tens\longrightarrow\Tens'$ be a natural mapping. It means that for any pair of objects $X,Y\in\cst_{G}$ we have morphism 
\[
\Phi^{XY}\in\Mor_{G}(X\Tens Y,X\Tens'Y)
\]
 and that for any pair of morphisms $r\in\Mor_{G}(X,Z)$ and $s\in\Mor_{G}(Y,T)$ the diagram
\[
\etyk{Db}
\vcenter{
\xymatrix{
X\Tens Y\ar[rr]^-{r\Tens s}\ar[d]_-{\Phi^{XY}}&&Z\Tens T\ar[d]^-{\Phi^{ZT}}
\\
X\Tens'Y\ar[rr]_-{r\Tens's}&&Z\Tens'T
}}
\]
is commutative. We know that $\zesp\Tens X=\zesp\Tens'X=X=X\Tens\zesp$ $=X\Tens'\zesp$. Therefore $\Phi^{X\zesp},\Phi^{\zesp X}\in\Mor_{G}(X,X)$. We say that $\Phi$ is normalized if $\Phi^{X\zesp}=\id_{X}=\Phi^{\zesp X}$ for any $X\in\cst_{G}$. In general setting $\varphi^{X}=\Phi^{X\zesp}$ and $\psi^{X}=\Phi^{\zesp X}$ we obtain two natural maps $\varphi$ and $\psi$ acting from $\id_{\cst_{G}}$ into itself. \vs
 
Let $\Tens,\Tens'$ be monoidal structures on $\cst_{G}$. We denote by $\alpha$, $\beta$, $\alpha'$ and  $\beta'$ the corresponding natural mappings:

\[
\begin{array}{r@{\;=\;}c@{\;}c@{\;}c@{\;\in\Mor_{G}(}c@{,X}c@{Y)}l}
\alpha^{XY}&\id_{X}&\Tens&1_{Y}&X&\Tens&,\\
\beta^{XY}&1_{X}&\Tens&\id_{Y}&Y&\Tens&,\\
\alpha'^{XY}&\id_{X}&\Tens'&1_{Y}&X&\Tens'&,\\
\beta'^{XY}&1_{X}&\Tens'&\id_{Y}&Y&\Tens'&
\end{array}
\]
for any $X,Y\in\cst_{G}$.\vs

 Let $\Phi:\Tens\longrightarrow\Tens'$ be a normalised natural mapping. Replacing in \rf{Db} $X,Y,Z,T,r,s$ by $X,\zesp,X,Y,\id_{X},1_{Y}$ respectively we obtain $\alpha'^{XY}=\Phi^{XY}\comp\alpha^{XY}$. Similarly replacing $X,Y,Z,T,r,s$ by $\zesp,Y,X,Y,1_{X},\id_{Y}$ we get $\beta'^{XY}=\Phi^{XY}\comp\beta^{XY}$. This way we showed that for any normalised natural mapping $\Phi:\Tens\longrightarrow\Tens'$ we have a commutative diagram
\[
\etyk{NNM}
\vcenter{
\xymatrix{
&\Tens\ar[dd]^-{\Phi}&\\ 
{\rm Proj}_1\ar[ru]^-{\alpha}\ar[rd]_-{\alpha'}&&{\rm Proj}_2\ar[lu]_-{\beta}\ar[ld]^-{\beta'}\\
&\Tens'&
}}.
\]

\begin{Thm}
\label{j}
Let $\Tens$ and $\Tens'$ be monoidal structures on $\cst_{G}$ corresponding to the same $R$-matrix.~Then there exists one and only one normalized natural mapping $\Phi:\Tens\longrightarrow\Tens'$.  For any $X,Y\in\cst_{G}$ the morphism 
\[
\Phi^{XY}\in\Mor_{G}(X\Tens Y,X\Tens' Y)
\]
 is an isomorphism.
\end{Thm}

 {\bf Plan of the proof of Theorem \ref{j}.}
Let $\Tens$ and $\Tens'$ be two monoidal structures \linebreak on $\cst_{G}$ with Properties 1, 2 and 3. Assume that the corresponding $R$-matrices coincide: $R=R'$. Our aim is to find a natural mapping $\Phi$ such that \rf{NNM}
is a commutative diagram. To show Theorem \ref{j} is sufficient to  prove the following two Propositions:

\begin{Prop}
\label{j1}
For any $X,Y\in\cst_{G}$ there exists unique isomorphism 
\[\Phi^{XY}\in\Mor(X\Tens Y,X\Tens' Y)
\]
 such that
\[
\etyk{D9}
\vcenter{
\xymatrix{
&X\Tens Y\ar[dd]^-{\Phi^{XY}}&\\ 
X\ar[ru]^-{\alpha^{XY}}\ar[rd]_-{\alpha'^{XY}}&&Y\ar[lu]_-{\beta^{XY}}\ar[ld]^-{\beta'^{XY}}\\
&X\Tens' Y&
}}
\]
is a commutative diagram. 
\end{Prop}

\begin{Prop}
\label{j2}
The collection of isomorphisms 
\[
\left(\Phi^{XY}\right)_{X,Y\in\cst_{G}}
\]
 introduced in Proposition \ref{j1} gives rise to a normalized natural mapping $\Phi:\Tens\rightarrow\Tens'$. \Vs{5}
\end{Prop}

\section{A model for $A\Tens A$}
\label{Sec10}

At first we shall prove Proposition \ref{j1} for $X=Y=A$. To this end we have to build a model for $A\Tens A$ independent of the choice of $\Tens$. \vs

Let $(\pi,\pihat)$ be a Heisenberg pair acting on a Hilbert space $H$ (see section \ref{Sec2}). Representations are morphisms into the algebra of all compact operators. Therefore $\pi\in\Mor(A,B_{0}(\Hil))$ and $\pihat\in\Mor(\Ahat,B_{0}(\Hil))$, where $B_{0}(\Hil)$ is the algebra of all compact operators acting on $\Hil$. We shall consider the following diagram in the category $\cst$:

\[
\etyk{D3}
\vcenter{
\xymatrix{
&A^{\Tens2}\ar[dd]^-{\skrD}&\\ 
& &\\
A\ar[r]^-{\skrR}\ar[ruu]^-{\alpha^{AA}}\ar[rdd]_-{\gamma}&\hspace{2mm}A^{\Tens2}\tens B_{0}(\Hil)\hspace{2mm}&A\ar[l]^-{\skrL}\ar[luu]_-{\beta^{AA}}\ar[ldd]^-{\delta}\\ 
&&\\
&A\tens B_{0}(\Hil)\ar[uu]^-{\skrU}&
}},
\]
where
\[
\etyk{wzorki}
\left\{
\begin{array}{c}
\begin{array}{c@{(a)=\;}cc@{(a)=\;}c@{\Delta_{R}(a)}l}
\skrR&\I_{A^{\Tens2}}\tens\,\pi(a),&
\skrL&(\beta^{AA}\tens\pihat)&,\\
\gamma&\I_{A}\tens\,\pi(a),&
\delta&(\id_{A}\tens\,\pihat)&,\Vs{5}
\end{array}\\
\skrU=\beta^{AA}\tens\id_{B_{0}(\Hil)}\Vs{5}
\end{array}
\right.
\]
and $\Delta_{R}$ is the morphism introduced by \rf{DeltaR}. Clearly $\skrU$ is injective. The description of the `down' morphism $\skrD: A^{\Tens2}\longrightarrow A^{\Tens2}\tens B_{0}(\Hil)$
 is more complicated. It is the composition of four injective morphisms (reading from right):
 \[
 \begin{array}{c@{\,:\,}r@{\;\longrightarrow\;}l}
 \Delta\Tens\id_{A}&A^{\Tens2}&(A\tens A)\Tens A=A\tens A^{\Tens2},\\ \Vs{5}
 \pi\tens\id_{A^{\Tens2}}&A\tens A^{\Tens2}&B_{0}(\Hil)\tens A^{\Tens2},\\ \Vs{5}
\flip&\ B_{0}(\Hil)\tens A^{\Tens2}&A^{\Tens2}\tens B_{0}(\Hil),\\ \Vs{5}
 \Ad_{\Vhat_{\alpha\pihat}}&A^{\Tens2}\tens B_{0}(\Hil)&A^{\Tens2}\tens B_{0}(\Hil),
\end{array}
\]
where $\Ad_{\Vhat_{\alpha\pihat}}(x)=\Vhat_{\alpha\pihat}\,x\,\Vhat_{\alpha\pihat}^{*}$ ($x\in A^{\Tens2}\tens B_{0}(\Hil)$, $\Vhat=\flip(\!V^{*}\!)$) and the equality of $\cst$-algebras in the first row comes from \rf{mixed}. So we have:
\[
\skrD= \Ad_{\Vhat_{\alpha\pihat}}\comp\flip\comp( \pi\tens\id_{A^{\Tens2}})\comp(\Delta\Tens\id_{A}).
\]
 
Lower triangles of \rf{D3} are obviously commutative. Upper triangles require a longer computation. Inserting in \rf{Nat} $X=Y=Y'=A$, $X'=A\tens A$, $\varphi=\Delta$ and $\psi=\id_{A}$  and using \rf{mi} we get
\[
\etyk{66}
\left\{
\begin{array}{r@{\;=\;}c@{\;=\;}c}
(\Delta\Tens\id_{A})\comp\alpha^{AA}&\alpha^{A\tens A,A}\comp\Delta&(\id_{A}\tens\,\alpha^{AA})\comp\Delta, \\
(\Delta\Tens\id_{A})\comp\beta^{AA}&\beta^{A\tens A,A}&1_{A}\tens\beta^{AA}.\Vs{5}
\end{array}
\right.
\]
 Let $\Vhat=\flip(V^{*})\in\M(A\tens\Ahat)$. Assigning the leg 1 to $A$, leg 2 to $A^{\Tens2}$ and leg 3 to $\Ahat$ we may rewrite \rf{5} in the form $V_{1\beta}\Rhat_{13}=
\Vhat_{\alpha3}V_{1\beta}\Vhat_{\alpha3}\Vs{4}^{*}$. Applying $\id_{\Ahat}\tens\beta^{AA}\tens\id_{\Ahat}$ to the both sides of \rf{DeltaR} we get
\[
%\etyk{Rhatt}
\left(\Vs{4}\id_{\Ahat}\tens\,(\beta^{AA}\tens\id)\Delta_{R}\right)V=V_{1\beta}\Rhat_{13}=\Vhat_{\alpha3}V_{1\beta}\Vhat_{\alpha3}\Vs{4}^{*}.
\]
Therefore for any $a\in A$ we have
\[
%\etyk{Rhat}
(\beta^{AA}\tens\id)\Delta_{R}(a)=\Vhat_{\alpha2}(\beta^{AA}(a)\tens I)\Vhat_{\alpha2}\Vs{4}^{*}.
\]
Applying to the both sides morphism $\id_{A\Tens A}\tens\,\pihat$ and using the notation introduced above (see \rf{wzorki}) we get
\[
\etyk{67}
\skrL(a)=\Ad_{\Vhat_{\alpha\pihat}}(\beta^{AA}(a)\tens\I_{\Hil}).
\]
On the other hand starting with the second formula of \rf{66} we compute:
\[
\begin{array}{r@{\;=\;}l}
(\Delta\Tens\id_{A})\comp\beta^{AA}(a)&\I_{A}\tens\,\beta^{AA}(a),\\
\Vs{5}(\pi\tens\id_{A^{\Tens2}})\comp(\Delta\Tens\id_{A})\comp\beta^{AA}(a)&\I_{\Hil}\tens\,\beta^{AA}(a),\\
\Vs{5}\flip\comp\,(\pi\tens\id_{A^{\Tens2}})\comp(\Delta\Tens\id_{A})\comp\beta^{AA}(a)&\beta^{AA}(a)\tens\I_{\Hil}.
\end{array}
\]
Inserting this result into \rf{67} we see that $\skrL(a)=\skrD\comp\beta^{AA}(a)$. It shows that the upper right triangle in \rf{D3} is commutative. To prove that upper left triangle is commutative we start with formula \rf{Delta}. Applying to the both sides $\id_{B_{0}(\Hil)}\tens\alpha^{AA}$ and using the first formula of \rf{66} we obtain
\[
\begin{array}{r@{\;=\;}c}
(\pi\tens\id_{A^{\Tens2}})\comp(\id_{A}\tens\alpha^{AA})\comp\Delta(a)&V_{\pihat\alpha}(\pi(a)\tens\I_{A^{\Tens2}})V_{\pihat\alpha}^{*},\\
\Vs{5}(\pi\tens\id_{A^{\Tens2}})\comp(\Delta\Tens\id_{A})\comp\alpha^{AA}(a)&V_{\pihat\alpha}(\pi(a)\tens\I_{A^{\Tens2}})V_{\pihat\alpha}^{*},\\
\Vs{5}\flip\comp(\pi\tens\id_{A^{\Tens2}})\comp(\Delta\Tens\id_{A})\comp\alpha^{AA}(a)&\Vhat_{\alpha\pihat}^{*}(\I_{A^{\Tens2}}\tens\,\pi(a))\Vhat_{\alpha\pihat},\\
\Vs{5}\Ad_{\Vhat_{\alpha\pihat}}\comp\flip\comp(\pi\tens\id_{A^{\Tens2}})\comp(\Delta\Tens\id_{A})\comp\alpha^{AA}(a)&\I_{A^{\Tens2}}\tens\,\pi(a),\\
\Vs{5}\skrD\comp\alpha^{AA}(a)&\skrR(a)
\end{array}
\]

\Vs{6}This way we proved that \rf{D3} is a commutative diagram. Proposition \ref{P1} shows now that there exists unique injective $\varphi\in\Mor(A^{\Tens2},A\tens B_{0}(\Hil))$ such that the diagram
\[
\vcenter{
\xymatrix{
&A\Tens A\ar[dd]^-{\varphi}&\\ 
& &\\
A\ar[r]_-{\gamma}\ar[ruu]^-{\alpha^{AA}}&\hspace{2mm}A\tens B_{0}(\Hil)\hspace{2mm}&A\ar[l]^-{\delta}\ar[luu]_-{\beta^{AA}}
}},
\]
is commutative. The reader should notice that the lower row of this diagram is independent of the monoidal structure $\Tens$. One could say that $(\gamma,\delta,\gamma(A)\delta(A))$ is a model for $(\alpha^{AA},\beta^{AA},A\Tens A)$.\vs

If $\Tens'$ is another monoidal structure on $\cst_{G}$ then combining the above diagram with the one for $\Tens'$ we obtain the commutative diagram
\[
\vcenter{
\xymatrix{
&A\Tens A\ar[dd]^-{\varphi}&\\ 
& &\\
A\ar[r]_-{\gamma}\ar[ruu]^-{\alpha^{AA}}\ar[rdd]_-{\alpha'^{AA}}&\hspace{2mm}A\tens B_{0}(\Hil)\hspace{2mm}&A\ar[l]^-{\delta}\ar[luu]_-{\beta^{AA}}\ar[ldd]^-{\beta'^{AA}}\\ 
&&\\
&A\Tens' A\ar[uu]_-{\varphi'}&
}},
\]

Using again Proposition \ref{P1} we find unique isomorphism $\Phi^{AA}\in\Mor(A\Tens A,A\Tens'A)$ such that the diagram
\[
\etyk{D8}
\vcenter{
\xymatrix{
&A\Tens A\ar[dd]^-{\Phi^{AA}}&\\ 
A\ar[ru]^-{\alpha^{AA}}\ar[rd]_-{\alpha'^{AA}}&&A\ar[lu]_-{\beta^{AA}}\ar[ld]^-{\beta'^{AA}}\\ 
&A\Tens' A&
}}.
\]
is commutative. This way we proved that Proposition \ref{j1} holds for $X=Y=A$.

\section{Proof of Theorem \ref{j}}
\label{Sec11}

\begin{pf}[ of Proposition \ref{j1}]
Let $X,Y\in\cst_{G}$. Then $\rho^{X}\Tens\rho^{Y}$ is a morphism from $X\Tens Y$ into $(X\tens A)\Tens(Y\tens A)$. The latter $\cst$-algebra equals to $X_{\tr}\Tens A\Tens Y_{\tr}\Tens A$. Applying $\flip$ to $A\Tens Y_{\tr}$ we obtain a morphism
\[
t=(\id_{X}\Tens\flip\Tens\id_{A})\comp(\rho^{X}\Tens\rho^{Y})
\]
acting from $X\Tens Y$ into $X\tens Y\tens A^{\Tens2}$. Let 
\[
\begin{array}{r@{\;=\;t\,\comp\,}c@{^{XY}\in\Mor(}c@{,X\tens Y\tens A^{\Tens2})}l}
r&\alpha&X&,\\
s&\beta&Y&.
\end{array}
\]
Then
\[
\begin{array}{r@{\;=\;}l}
r&(\id_{X}\Tens\flip\Tens\id_{A})\comp(\rho^{X}\Tens\rho^{Y})\comp\alpha^{XY}\\
&(\id_{X}\Tens\flip\Tens\id_{A})\comp\alpha^{X\tens A,Y\tens A}\comp\rho^{X}\\
&(\id_{X}\tens 1_{Y}\tens\alpha^{AA})\comp\rho^{X},
\end{array}
\]
\[
\begin{array}{r@{\;=\;}l}
s&(\id_{X}\Tens\flip\Tens\id_{A})\comp(\rho^{X}\Tens\rho^{Y})\comp\beta^{XY}\\
&(\id_{X}\Tens\flip\Tens\id_{A})\comp\beta^{X\tens A,Y\tens A}\comp\rho^{Y}\\
&(1_{X}\tens\id_{Y}\tens\beta^{AA})\comp\rho^{Y},
\end{array}
\]

Replacing in the above formulae $\Tens$ be $\Tens'$ we obtain morphisms $r'$, $t'$ and $s'$ acting from $X$, $X\Tens'Y$ and $Y$ into $X\tens Y\tens(A^{\Tens'2})$ such that 
\[
\begin{array}{r@{\;=\;}l}
t'&(\id_{X}\Tens'\flip\Tens'\id_{A})\comp(\rho^{X}\Tens'\rho^{Y}),\\
r'&t'\comp\alpha'^{XY}=(\id_{X}\tens 1_{Y}\tens\alpha'^{AA})\comp\rho^{X},\\
s'&t'\comp\beta'^{XY}=(1_{X}\tens\id_{Y}\tens\beta'^{AA})\comp\rho^{Y},
\end{array}
\]

\Vs{10}This way we constructed a commutative diagram:

\[
\etyk{ok}
\vcenter{
\xymatrix{
&&X\Tens Y\ar[dd]^-{t}&&\\
&&&&\\
&&X\tens Y\tens A^{\Tens2}\ar[dd]^-{\id_{_{X\tens Y}}\!\tens\Phi^{AA}}&&\\
X\ar@/^3pc/[uuurr]^-{\alpha^{XY}}
\ar[urr]^-{r}\ar[drr]_-{r'}
\ar@/_3pc/[dddrr]_-{\alpha'^{XY}}&&&&
Y\ar@/_3pc/[uuull]_-{\beta^{XY}}
\ar[ull]_-{s}\ar[dll]^-{s'}
\ar@/^3pc/[dddll]^-{\beta'^{XY}}\\
&&X\tens Y\tens A^{\Tens'2}&&\\
&&&&\\
&&X\Tens'Y\ar[uu]_-{t'}&&
}}
\]
\vs

 Commutativity of the triangles with straight arrow sides immediately follows from diagram \rf{D8}\vs

Removing from \rf{ok} $r$, $X\tens Y\tens A^{\Tens2}$ and $s$ and replacing $\id_{_{X\tens Y}}\tens\Phi^{AA}$ and $t$ by their composition we obtain a commutative diagram of the form \rf{D1}. Proposition \ref{P1} shows now that there exists an isomorphism $\Phi^{XY}\in\Mor(X\Tens Y,X\Tens'Y)$ that makes diagram \rf{D9} commutative. Proposition \ref{j1} is proven in full generality.
\end{pf}\vs
 
\begin{pf}[ of Proposition \ref{j2}]
First we notice that $\Phi$ is normalised. Indeed inserting  in \rf{D9} $Y=\zesp$ we see that $\Phi^{X\zesp}=\id_{X}$. Similarly puting $X=\zesp$ we get $\Phi^{\zesp Y}=\id_{Y}$. \vs

Next we have to show that $\Phi^{XY}$ are morphisms in category $\cst_{G}$. It means that the diagram
\[
\etyk{Da}
\vcenter{
\xymatrix{
X\Tens Y\ar[rr]^-{\rho^{X\Tens Y}}\ar[d]_-{\Phi^{XY}}&&(X\Tens Y)\tens A\ar[d]^-{\Phi^{XY}\tens\id_{A}}
\\
X\Tens'Y\ar[rr]_-{\rho^{X\Tens'Y}}&&(X\Tens'Y)\tens A
}}
\]
is commutative for all $X,Y\in\cst_{G}$. We have:
\[
\begin{array}{r@{\;=\;}l}
(\Phi^{XY}\tens\id_{A})\comp\rho^{X\Tens Y}\comp\alpha^{XY}&
(\Phi^{XY}\tens\id_{A})\comp(\alpha^{XY}\tens\id_{A})\comp\rho^{X}\\
&(\Phi^{XY}\comp\alpha^{XY}\tens\id_{A})\comp\rho^{X}\\
&(\alpha'^{XY}\tens\id_{A})\comp\rho^{X}\\
&\rho^{X\Tens'Y}\comp\alpha'^{XY}\\
&\rho^{X\Tens'Y}\comp\Phi^{XY}\comp\alpha^{XY}
\end{array}
\]
\[
\begin{array}{r@{\;=\;}l}
(\Phi^{XY}\tens\id_{A})\comp\rho^{X\Tens Y}\comp\beta^{XY}&
(\Phi^{XY}\tens\id_{A})\comp(\beta^{XY}\tens\id_{A})\comp\rho^{Y}\\
&(\Phi^{XY}\comp\beta^{XY}\tens\id_{A})\comp\rho^{Y}\\
&(\beta'^{XY}\tens\id_{A})\comp\rho^{Y}\\
&\rho^{X\Tens'Y}\comp\beta'^{XY}\\
&\rho^{X\Tens'Y}\comp\Phi^{XY}\comp\beta^{XY}
\end{array}
\]\vs

\noindent Remembering that $X\Tens Y=\alpha^{XY}(X)\beta^{XY}(Y)$ we see that \rf{Da} is a commutative diagram.\vs
 
Finally we have to show that $\Phi$ is a natural mapping. Let $r\in\Mor_{G}(X,Z)$ and $s\in\Mor_{G}(Y,T)$. Then
\[
\begin{array}{r@{\;=\;}l}
\Phi^{ZT}\comp(r\Tens s)\comp\alpha^{XY}&\Phi^{ZT}\comp\alpha^{ZT}\comp r\\
&\alpha'^{ZT}\comp r\\
&(r\Tens's)\comp\alpha'^{XY}\\
&(r\Tens's)\comp\Phi^{XY}\comp\alpha^{XY}
\end{array}
\]
\[
\begin{array}{r@{\;=\;}l}
\Phi^{ZT}\comp(r\Tens s)\comp\beta^{XY}&\Phi^{ZT}\comp\beta^{ZT}\comp s\\
&\beta'^{ZT}\comp s\\
&(r\Tens's)\comp\beta'^{XY}\\
&(r\Tens's)\comp\Phi^{XY}\comp\beta^{XY}
\end{array}
\]
Remembering that $X\Tens Y=\alpha^{XY}(X)\beta^{XY}(Y)$ we see that \rf{Db} is a commutative diagram. $\Phi$ is a natural mapping.  Proposition \ref{j2} is shown.
\end{pf}\vs
 
This way we proved Theorem \ref{j}.

\section{A remark on Property 3.}
\label{Sec12}

In Theorem \ref{MRW} Properties 1, 2 and 3 appears as the part of the statement. Therefore they are formulated in the strongest version. On the other hand in Theorem \ref{S} they belong to the assumptions and it is desirable to formulate them in a possibly weak form.\vs

We shall use the two-dimensional abelian $\cst$-algebra $D=\zesp^{2}$ with the trivial action of $G$. Then $D\in\cst_{G}$ and
\[
\rho^{D}(r)=r\tens\I_{A}
\]
for any $r\in D$. Let $p\in D$ be one of the two nontrivial projections in $D$. Then $q=I_{D}-p$ is the second nontrivial projection and $D$ is the linear span of $\left\{p,q\right\}$. It turns out that Property 3 may be replaced by apparently weaker\footnote{This formulation of Property 3 belongs to Ralf Meyer.}\vs

{\bf Property 3':}
$\Tens$ reduces to $\tens$, when one of the involved $\cst$-algebras is $D$. More precisely: For any 
$X\in\cst_{G}$ we have: $X\Tens D=X\tens D$ as $\cst$-algebras and
\[
\begin{array}{r@{\;=\;}l}
\alpha^{XD}(x)&x\tens I_{D},\\
\beta^{XY}(r)&I_{X}\tens r
\end{array}
\]
for any $x\in X$ and $r\in D$. Similarly for any 
$Y\in\cst_{G}$ we have: $D\Tens Y=D\tens Y$ as $\cst$-algebras and
\[
\begin{array}{r@{\;=\;}l}
\alpha^{DY}(r)&r\tens I_{Y},\\
\beta^{DY}(y)&I_{D}\tens y
\end{array}
\]
for any $r\in D$ and $y\in Y$.\vs

One may also use the statement of the Proposition \ref{cp} as possible replacement:\vs

{\bf Property 3'':} For any $X,Y\in\cst_{G}$ and any $x\in X$ and $y\in Y$, elements $\alpha^{XY}(x)$ and $\beta_{XY}(y)$ commute if one of the elements $x$, $y$ is $G$-invariant.

\begin{Thm}
%\label{}
Let $G=(A,\Delta)$ be a locally compact quantum group and $\Tens$ be a monoidal structure on $\cst_{G}$ having Properties 1 and 2. Then Properties 3, 3' and 3'' are equivalent.\vs
\end{Thm}
\begin{pf}
It is obvious that $3\Rightarrow 3'$.  We also know (cf. Proposition \ref{cp}) that $3\Rightarrow 3''$. We shall prove the converse implications. Let $X,Y\in\cst_{G}$ and the action $\rho^{X}$ be trivial. We may assume that $X\subset B(K)$, where $K$ is a Hilbert space. Denote the algebra of all compact operators acting on $K$ by $B_{0}(K)$ and provide it with the trivial action of $G$. Then $B_{0}(K)\in\cst_{G}$ and the embedding $i:X\hookrightarrow B(K)$ is a $\cst_{G}$ morphism from $X$ into $B_{0}(K)$: $i\in\Mor_{G}(X,B_{0}(K))$. Therefore $i\Tens\id_{Y}\in\Mor_{G}(X\Tens Y,B_{0}(K)\Tens Y)$ is an injective morphism. We shall identify $X\Tens Y$ with its image:
 \[
 X\Tens Y\subset \M(B_{0}(K)\Tens Y).
 \]
 Then in virtue of diagram \rf{Nat} (with $X',Y',\varphi,\psi$ replaced by $B_{0}(K),Y,i,\id_{Y}$) we have
 \[
 \begin{array}{r@{\;=\;}l}
\alpha^{XY}(x)&\alpha^{B_{0}(K),Y}(x),\\
 \beta^{XY}(y)&\beta^{B_{0}(K),Y}(y)
\end{array}
 \]
for any $x\in X$ and $y\in Y$. \vs

 Choosing a faithful representation we may assume that  the $\cst$-algebra $B_{0}(K)\Tens Y$ is contained in $B(H)$, where  $H$ is a Hilbert space. Then $\alpha^{B_{0}(K),Y}$ is a representation of $B_{0}(K)$ acting on $H$.\vs 
   
Using the well known property of the algebra of all compact operators (\cite{Arv}, Corollary 1, page 20) we see that $H$ is of the form $H=K\tens H'$ ($H'$ is another Hilbert space) and
\[
\etyk{2053}
\alpha^{B_{0}(K),Y}(x)=x\tens\I
\] 
for any $x\in B_{0}(K)$.\vs

Let us fix a nontrivial projection $p\in D$. Any orthonormal projection $x\in B_{0}(K)$ is of the form  $x=\varphi(p)$, where $\varphi\in\Mor_{G}(D, B_{0}(K))$. By Property 3', $\alpha^{DY}(p)$ commutes with $\beta^{DY}(y)$. Therefore $\left(\varphi\Tens\id_{Y}\right)\alpha^{DY}(p)=\alpha^{B_{0}(K),Y}(x)$ commutes with $\left(\varphi\Tens\id_{Y}\right)\beta^{DY}(y)=\beta^{B_{0}(K),Y}(y)$. Remembering that the algebra $B_{0}(K)$ coincides with the closed linear span of all its orthogonal projections we conclude that the commutator
\[
\etyk{2052}
\left[\alpha^{B_{0}(K),Y}(x),\beta^{B_{0}(K),Y}(y)\right]=0
\]
for any $x\in B_{0}(K)$ and $y\in Y$. The reader should notice that using Property 3'' (instead of 3') one obtains the same result.\vs

Formulae \rf{2053} and \rf{2052} show that $\beta^{B_{0}(K),Y}(y)\in\I\tens B(H')$. Therefore there exists a faithful representation $\pi\in\Rep(Y,H')$ such that $\beta^{B_{0}(K),Y}(y)=\I\tens\pi(y)$. Identifying $Y$ with $\pi(Y)$ we have
\[
\beta^{B_{0}(K),Y}(y)=\I\tens\, y
\]
for all $y\in Y$. Combining formulae obtained so far we get
\[
 \begin{array}{r@{\;=\;}l}
\alpha^{XY}(x)&x\tens\I,\\
 \beta^{XY}(y)&\I\tens\, y
\end{array}
 \]
for any $x\in X$ and $y\in Y$. Now formula \rf{wtorek} shows that $X\Tens Y=X\tens Y$. The case, when the action $\rho^{Y}$ is trivial may be treated in the same way. The implications $3'\Rightarrow 3$ and $3''\Rightarrow 3$ are proved. 
\end{pf}

\section{Category $_{G}\cst$}
\label{Sec13}

This is the category of $\cst$-algebras equipped with left actions of a locally compact quantum group $G=(A,\Delta)$. Let $X$ be a $\cst$-algebra and $\lambda\in\Mor(X,X\tens A)$. We say that $\lambda$ is a left action of $G$ on $X$ if

1. The diagram
\[
\etyk{la}
\vcenter{
\xymatrix{
X\ar[rr]^-{\lambda}\ar[d]_-{\lambda}&&A\tens X\ar[d]^-{\id\tens\lambda}\\
A\tens X\ar[rr]_-{\Delta\tens\id}&&A\tens A\tens X
}}
\]
is a commutative,\vs

2. $\ker\lambda=\zero$,\vs

3. $(A\tens\I)\lambda(X)=A\tens X$ (Podle\'s condition).\vs

Let $X,Y$ be $\cst$-algebras equipped with left actions $\lambda^{X},\lambda^{Y}$ of $G$ and $\gamma\in\Mor(X,Y)$. We say that $\gamma$ intertwines the actions of $G$ if the diagram
\[
\etyk{MorGL}
\vcenter{
\xymatrix{
X\ar[rr]^{\gamma}\ar[d]_-{\lambda^{X}}&&Y\ar[d]^-{\lambda^{Y}}
\\
A\tens X\ar[rr]_-{\id\tens\gamma}&&A\tens Y
}}
\]
is commutative. By definition morpisms in the category $_{G}\cst$ are morphisms in the category $\cst$ intertwining the actions of $G$. The set of all morphisms from $X$ to $Y$ will be denoted by $_{G}\Mor(X,Y)$.\vs

Let $\Tens$ be a monoidal structure on $_{G}\cst$. As in the section \ref{Sec5} one can use formulae \rf{dfg} to introduce natural mappings $\alpha$ and $\beta$ acting from ${\rm Proj}_1$ and ${\rm Proj}_2$ into $\Tens$. Then  for any $X,Y\in\,_{G}\cst$ we have morphisms $\alpha^{XY}\in\,_{G}\Mor(X,X\Tens Y)$ and $\beta^{XY}\in\,_{G}\Mor(Y,X\Tens Y)$. In this context we shall also use the abbreviation \rf{abr}. 

\begin{Thm}
\label{MRWL}
Let $G=(A,\Delta)$ be a quasitriangular locally compact quantum group with a unitary $R$-matrix $R\in \M(\Ahat\tens\Ahat)$. Then there exists  a monoidal structure $\Tens$ on $_{G}\cst$ having Properties 1, 2 and 3 and such that
\[
\etyk{RaH}
V_{2\beta}V_{1\alpha}=R_{12}V_{1\alpha}V_{2\beta}
\]
\end{Thm}

\begin{Thm}
\label{SL}
Let $G=(A,\Delta)$ be a locally compact quantum group and $\Tens$ be a monoidal structure on $_{G}\cst$ having Properties 1, 2 and 3. Then there exists $($unique$)$ unitary $R$-matrix $R\in\M(\Ahat\tens\Ahat)$ such that
\[
V_{2\beta}V_{1\alpha}=R_{12}V_{1\alpha}V_{2\beta}
\]
\end{Thm}

\begin{pf}
It is easy to reduce these theorems to Theorems \ref{MRW} and \ref{S}. To this end we consider the group $G^{\text opp}=(A,\flip\comp\Delta)$ opposite to $G$ (see section \ref{Sec2}). Clearly any left action $\lambda$ of $G$ on a $\cst$-algebra defines the corresponding right action $\rho=\flip\comp\lambda$ of $G^{\rm opp}$ on the same algebra and vice-versa. It means that the categories $_{G}\cst$ and $\cst_{G^{\rm opp}}$ coincides. In particular they have the same monoidal structures.\vs

We know (see the end of section \ref{Sec2}) that the passage from $G$ to $G^{\rm opp}$ consists in replacing $V$ and $R$ by $V^{*}$ and $R^{*}$. To end the proof the reader should notice that replacing in \rf{RH} $V$ and $R$ by $V^{*}$ and $R^{*}$ (and taking the $^{*}$ of both sides) we obtain \rf{RaH}.
\end{pf}\vs

By the same argument we may replace $\cst_{G}$ by $_{G}\cst$ in Theorem \ref{j}. It shows that the monoidal structure on $_{G}\cst$ is uniquely determined by $R$-matrix.
 
\section*{Acknowlegdement}
The author is very grateful to Ralf Meyer and Sutanu Roy for interesting and sti\-mulating discussions.

%\bibliographystyle{plain}
%\bibliography{SLWbib,Spis}

\end{document}